\documentclass[a4paper]{article}

\usepackage{amsfonts,amssymb,amsthm}
\usepackage[totalwidth=17cm,totalheight=22cm]{geometry}

\newtheorem{prop}{Proposition}[section]
\newtheorem{df}[prop]{Definition}
\newtheorem{lemma}[prop]{Lemma}
\newtheorem{thm}[prop]{Theorem}
\newtheorem{cor}[prop]{Corollary}

\newtheorem{remark}[prop]{Remark}

\newcommand{\II}{I\hspace{-0.1cm}I}
\newcommand{\III}{I\hspace{-0.1cm}I\hspace{-0.1cm}I}
\newcommand{\dr}{\partial}
\newcommand{\db}{\overline{\partial}}

\newcommand{\tr}{\mbox{tr}}
\newcommand{\be}{\begin{eqnarray}}
\newcommand{\ee}{\end{eqnarray}}
\newcommand{\C}{{\mathbb C}}
\newcommand{\R}{{\mathbb R}}
\newcommand{\cS}{{\mathcal S}}
\newcommand{\cT}{{\mathcal T}}

\begin{document}

\title{On the renormalized volume of hyperbolic 3-manifolds}
\date{November 2006 (v3)}
\author{Kirill Krasnov\thanks{School of Mathematical Sciences, University of
    Nottingham, Nottingham, NG7 2RD, UK} 
~and Jean-Marc Schlenker\thanks{
Laboratoire Emile Picard, UMR CNRS 5580,
Institut de Math{\'e}matiques, 
Universit{\'e} Paul Sabatier,
31062 Toulouse Cedex 9,
France.
\texttt{http://www.picard.ups-tlse.fr/\~{ }schlenker}.}
}
\maketitle

\begin{abstract} The renormalized volume of hyperbolic manifolds is a quantity 
motivated by the AdS/CFT correspondence of string theory and computed 
via a certain regularization procedure. The main aim of the present paper is
to elucidate its geometrical meaning. We use another regularization procedure 
based on surfaces equidistant to a given convex surface $\partial N$. The
renormalized volume computed via this procedure is equal 
to what we call the $W$-volume of the convex region $N$
given by the usual volume of $N$ minus the quarter of the integral of the mean
curvature over $\partial N$. The $W$-volume satisfies some remarkable
properties.  
First, this quantity is self-dual in the sense explained in the paper.
Second, it verifies some simple variational formulas analogous to the classical
geometrical Schl\"afli identities. These variational formulas are invariant
under a certain  
transformation that replaces the data at $\partial N$ by those at infinity of
$M$.  
We use the variational formulas in terms of the data at infinity to give
a simple geometrical proof of results of Takhtajan et al on the K\"ahler
potential on various moduli spaces.
\end{abstract}

\section{Introduction}

\paragraph{The renormalized volume.}

In this paper we study the so-called renormalized volume of hyperbolic
3-manifolds whose definition is 
motivated by the AdS/CFT correspondence of string theory
\cite{Witten:1998qj}. In the context of this                        
correspondence one is interested in computing the gravity action 
\be\label{gr-action}
S_{gr}[g] = \frac{1}{2} \int_M (R-2\Lambda)dv + \int_{\partial M} Hda
\ee
for an Einstein metric $g$ on d-dimensional manifold $M$. Here we have set the
dimensionful Newton constant 
typically present in the action (\ref{gr-action}) to $8\pi G=1$, and $\Lambda$
is the cosmological constant, 
which is assumed to be negative. The quantity $R$ is the curvature scalar, $H$
is the mean  
curvature of the boundary $\partial M$, and $dv, da$ are the volume and area
forms correspondingly. 

In the context of AdS/CFT correspondence the manifold $M$ is non-compact with a conformal boundary $\partial M$. 
The metric $g$ on $M$ determines the conformal class of the boundary $\partial M$. 
One would like to compute $S_{gr}$ as a functional
of this conformal class. However, because the manifold
$M$ is non-compact the functional $S_{gr}[g]$ diverges. One notices, however, that
this divergence is of a special type. Thus, let us introduce a compact sub-manifold $N\subset M$
and compute the action $S_{gr}[g,N]$ (\ref{gr-action}) inside $N$. As one sends $\partial N$ towards
the infinity of $M$ one can show the divergent part of $S_{gr}[g,N]$ is given by an integral
over $\partial N$ of a local quantity expressible in terms of the first and second fundamental forms
of $\partial N$. This suggests an idea of {\it renormalization} in which these divergent quantities
are subtracted, after which the limit $\partial N\to \partial M$ can be taken. This idea works 
in any dimension, see \cite{Balasubramanian:1999re} for the required subtraction procedure. 

\paragraph{3-dimensional manifolds.}

In this paper we are interested in 
the simplest case of 3-dimensional spaces. The simplifications arising in this case are
as follows. First, Einstein equations in 3 dimensions imply the metric $g$ to be of constant curvature.
Thus, in our case of a negative cosmological constant we are led to consider constant negative curvature
manifolds. In 3 dimensions the radius of curvature $l$ is related to the cosmological constant
in a simple way $l=1/\sqrt{-\Lambda}$. We shall set the radius of curvature $l=1$ in what follows, 
which in 3 dimensions is equivalent to $\Lambda=-1$. Thus, the 
Riemannian manifolds $(M, g)$ we are going
to consider are hyperbolic. As is well known from Bers' work \cite{Bers} on simultaneous uniformization and from
its generalization to arbitrary Kleinian manifolds, the manifolds $(M,g)$ are completely
characterized by the conformal structures of all the boundary components of $M$. We would like
to compute the gravity action $S_{gr}$ as a functional of the conformal structure of $\partial M$.
The constant curvature condition gives $R=-6$. Therefore, the gravity action
reduces to  
\be\label{gr-1}
S_{gr} = - 2 \int_M dv + \int_{\partial M} Hda.
\ee
It is now easy to show that the volume $V(N), N\subset M$ diverges as $(1/2) A(\partial N)$, where
$A(\partial N)$ is the area of the boundary of $N$. There is also the so-called logarithmic divergence
to be discussed below. The integrated mean
curvature diverges as $\int_{\partial N} Hda \sim 2 A(\partial N)$. Altogether we see that
the gravity action functional diverges as $A(\partial N)$. It seems natural, therefore, 
to introduce a {\it renormalized} gravity action given by
\be 
\label{gr-action-reg}
S_{gr}[g] = \frac{1}{2} \int_M (R-2\Lambda)dv + \int_{\partial M} (H-1)da.
\ee
This action can then be computed for any compact sub-manifold $N\subset M$. However, 
because of the logarithmic divergence present in the volume, the
limit $N\to M$ does not exist. It is this logarithmic divergence that causes all the
difficulties that we are now to discuss. 

In an even number of dimensions, when a more involved
but similar in principle subtraction procedure is used, there is no logarithmic divergence and
the limit can indeed be shown to exist and to be independent 
on how exactly the surfaces $\partial N$ are taken
to approach the conformal boundary $\partial M$. However, when the dimension
is odd, as is in the case of interest for us, the limit of (\ref{gr-action-reg}) does
not exist, except in the special situation when all the boundary components are of
genus one. In general, the volume of $N$ grows as the logarithm of the area of $\partial N$
times the Euler characteristic of $\partial N$. As we shall see below, one can
subtract this divergence as well, but the resulting renormalized gravity action (or volume) then
turns out to depend on the limiting procedure. This is why in odd dimensions the renormalized
gravity action fails to be a true invariant of $M$. It is precisely for this
reason that the concept of renormalized volume has not been developed by
the geometry community. Indeed, the fact that the volume of a hyperbolic manifold
grows as half of the area of its boundary has been known for a very long time and
is mentioned, for example, in Thurston's famous notes on the subject. It was also
realized however that the limit $V(N)-(1/2)A(\partial N)$ does not exist. The
geometry community instead concentrated on e.g. a canonically defined volume
of the convex core of $M$. The logarithmic divergence (referred to as the conformal
anomaly in the physics literature) is studied in more details in \cite{Skenderis}.

In \cite{Krasnov:2000zq} one of us studied the renormalized volume of Schottky
3-manifolds. It was shown that there is a way to define the limit of the
quantity $V(N)-(1/2)A(\partial N)$ by choosing a foliation of $M$ near its
boundary by a family of surfaces $S_\rho[\phi]$ parameterized by a 
``Liouville'' (real-valued) field on the conformal boundary of $M$. Once such
a family of surfaces is used, the limiting procedure becomes well-defined. Indeed,
one takes the quantity $V_\rho - (1/2)A_\rho$ and subtracts from it $2\pi\rho(g-1)$,
where $g$ is the genus of the boundary of $M$. It is then easy to see that the limit exists and defines
the renormalized volume $V_R(M,\phi)$, which in addition to being a function
of the conformal structure of $\partial M$ also depends on the Liouville
field $\phi$ on $\partial M$. It was moreover shown by an explicit
computation that $V_R(M,\phi)$ is equal to the Liouville action $S_L(\partial M, \phi)$
of Takhtajan and Zograf \cite{TZ-Schottky}. To get a ``canonical'' quantity that depends
only on the conformal structure of $\partial M$ one can evaluate $V_R(M,\phi)=S_L(\partial M, \phi)$
on the canonical Liouville field corresponding to the metric of constant curvature $-1$
in the conformal class of $\partial M$. This canonical $\phi$ can also be obtained by
extremizing the functional $V_R(M,\phi)$ keeping the area of $\partial M$ as determined
by $\phi$ fixed. This variational principle leads to the unique canonically defined $\phi$. 
One gets the ``canonical'' renormalized volume that depends only on the conformal
structure of $\partial M$.

This way of defining the limit can be generalized to an arbitrary Kleinian manifold.
This has been done in \cite{Takhtajan:2002cc}, where the limiting procedure
of \cite{Krasnov:2000zq} was also improved in the sense that an invariant
family of Epstein surfaces \cite{Eps} was used for regularization, see also an earlier paper by one of
the authors \cite{Krasnov:2001cu}, where the same Epstein family of surfaces was
used for regularization, but of the CS formulation of gravity instead. In all cases, the renormalized volume
obtained via the limiting procedure was shown to agree with the Liouville action on $\partial M$.
Variation of the renormalized volume = Liouville action under the
changes of the conformal structure of $\partial M$ were studied in \cite{TZ-Schottky}, 
\cite{Takhtajan:2002cc}. It was shown that in all the cases the renormalized volume 
is equal to the K\"ahler potential for the Weil-Petersson metric on the moduli
space of $\partial M$. For quasi-Fuchsian spaces this implies, e.g. the quasi-Fuchsian
reciprocity of McMullen \cite{McMullen}.

\paragraph{The renormalized volume and equidistant foliations.}

In the present paper we undertake a further study of the renormalized volume of
hyperbolic 3-manifolds. We show that the limiting procedure via which the
volume is defined can be somewhat de-mystified by considering for
regularization a family of surfaces equidistant to a given one, following an
idea already used by C. Epstein \cite{epstein-duke} (and more recently 
put to use in \cite{minsurf}). Thus, the main
idea of the present work is to obtain the renormalized volume by taking
a convex domain $N\subset M$, and compute the renormalized volume of $M$ {\it
  with respect to $N$} as 
\be\label{rv}
V_R(M,N)=V(N) + \lim_{\rho\to\infty} \left( V(\partial N, \partial N_\rho) - (1/2) A(\partial N_\rho) - \sum_i 2\pi \rho(g_i-1) \right),
\ee
where $V(\partial N, \partial N_\rho)$ is the volume between the boundary $\partial N$
of the domain $N$ and the surface $\partial N_\rho$ located a distance $\rho$ from $\partial N$.
The quantity $A(\partial N_\rho)$ is the area of the surface $\partial N_\rho$, the sum
in the last term is taken over all boundary components of $M$ and $g_i$ are the genera of
these boundary components. The convexity of
the domain $N$ ensures that the equidistant surfaces $\partial N_\rho$ exist all the way to
infinity. This ensures that the limit $\rho\to\infty$ can be taken. Similarly, using the
combination (\ref{gr-1}) and subtracting the area of $\partial N_\rho$, as well as the
term linear in $\rho$, one defines the renormalized gravity action. In both cases the 
limit exists and can be computed in terms of the volume or the gravity action for $N$, 
see below for the corresponding expressions.

The limiting procedure used in \cite{Krasnov:2001cu}, \cite{Takhtajan:2002cc} is an example
of the limiting procedure described above, for the Epstein surfaces \cite{Eps} are equidistant. Thus,
the renormalized volume of references \cite{Krasnov:2000zq}, \cite{Takhtajan:2002cc}
is an example of the renormalized volume (\ref{rv}) where $N$ is a compact domain of $M$
contained inside a particular Epstein surface. However, it is obvious that the
renormalized volume (\ref{rv}) is more general as the domain $N$ in (\ref{rv}) can
be an arbitrary convex domain. It is also clear that the renormalized volume
defined via (\ref{rv}) is ``the most general one''. Indeed, the only constraint that
enters into the definition (\ref{rv}) is that equidistant surfaces are used for
regularization. However, this seems to be a necessary requirement to be able 
to subtract the logarithmic divergence. Thus, there does not seem any other way
to define the renormalized volume. 

\paragraph{Renormalized volume as the $W$-volume.}

The starting point is a simple formula for the renormalized volume 
(\ref{rv}):
\be
V_R(M, N) = W(N) - \sum_i \pi(g_i-1)~, 
\ee
where the sum in the last term is taken over all the boundary components. Here
the W-volume is defined as
\be
W(N) := V(N) - \frac{1}{4} \int_{\dr N} H da.
\ee
This formula for $V_R$ is a special case of a formula found by C. Epstein 
\cite{epstein-duke} for the renormalized volume of hyperbolic manifolds 
in any dimension. 

Thus, the renormalized volume of $M$ with respect to $N$ is, apart from an
uninteresting term given by a multiple of the Euler characteristic of the
boundary, just the W-volume of the domain $N$. All our other results
concern this W-volume. Note already that $W(N)$ is not equal to the 
Hilbert-Einstein functional of $N$ with its usual boundary term, it
differs from it in the coefficient of the boundary term.

\paragraph{Variation formula.} We prove a formula
for the first variation of W under changes of the metric inside $N$.
This formula is a simple consequence of the formula obtained
by Igor Rivin and one of us in \cite{sem-era}. It reads
\be\label{intr-1}
\delta W(N) = \frac{1}{4} \int_{\partial N} \langle\delta \II - \frac{H}{2} \delta I, I\rangle  da.
\ee
Thus, this formula suggest that the $W$-volume of a domain $N$ is a complicated
functional of the shape of this domain, as well as of the manifold $M$. However,
as we shall show, this functional depends on all these data in a very specific
way, through a certain combination that we introduce below and refer to as
the metric ``at infinity''.

\paragraph{Self-duality.}
One of the most interesting properties of the W-volume is that it
is self-dual. Thus, we recall that the Einstein-Hilbert functional
\be\label{gr-2}
I_{EH}(N) := V(N) - \frac{1}{2} \int_{\dr N} H da
\ee
for a compact domain $N\subset H^3$ of the hyperbolic space
(note a different numerical factor in front of the second term) is nothing but the
dual volume. Thus, recall that there is a duality between objects in $H^3$ and
objects in $dS_3$, the 2+1 dimensional de Sitter space. Under this duality
geodesic planes in $H^3$ are dual to points in $dS_3$, etc. This duality between
domains in the two spaces is easiest to visualize for convex polyhedra
(see \cite{RH}), but the duality works
for general domains as well. The fact that (\ref{gr-2}) is the volume of the
dual domain 
then is a simple consequence of the Schl\"afli formula, analogous to
(\ref{intr-1}), see 
the main text below. 

Thus, we can write:
\be
{}^*V(N) = V({}^*N) = V(N) - \frac{1}{2} \int_{\dr N} H da
\ee
for the volume of the dual domain. This immediately shows that
\be
W(N) = V(N) - \frac{1}{4} \int_{\dr N} H da = \frac{V(N)+{}^*V(N)}{2}.
\ee
Thus, the $W$-volume is self-dual in that this quantity for $N$ is equal to this quantity
for the dual domain ${}^*N$: $W(N)=W({}^*N)$. In the main text we shall also verify the self-duality more
directly by applying the Legendre transform to $W(N)$.

\paragraph{The $W$-volume and the Chern-Simons formulation.} An interesting remark is that
there is a very simple expression for the $W$-volume in terms of the so-called
Chern-Simons formulation of 2+1 gravity \cite{Witten:2+1}. In this formulation
the gravity action in the so-called first order formalism, in which the independent
variables are the triads and the spin connection, is shown to be given by the
difference of two Chern-Simons actions. This is easily shown for the ``bulk'' term
of the gravity action, while to get the boundary term as in (\ref{gr-action}) one needs 
in addition a certain set of boundary terms in terms of the Chern-Simons connections.
Remarkably, the combination that plays this role in the $W$-volume, namely 
\be
\frac{1}{2} \int_M (R-2\Lambda) + \frac{1}{2} \int_{\partial M} K,
\ee
which is different from the usual Einstein-Hilbert action, is exactly what
appears naturally in the Chern-Simons formulation, without the need for
any boundary terms. We refer the reader to e.g. \cite{Krasnov:2001cu}
for a demonstration of this fact, see formula (3.7) of this reference
as well as the related discussion. It would be of interest to understand
the relation, if any, between the self-duality of the $W$-volume
and the fact that it has such a simple expression in the Chern-Simons
formulation.

\paragraph{Description ``from infinity''.}
The other important property of the W-volume is that it can be interpreted as a quantity
that depends on the metric ``at infinity'' only, instead of being a functional of the ``shape'' of 
the boundary $\partial N$ of $N\subset M$. To demonstrate this one introduces the first and
second fundamental forms $I^*, \II^*$ obtained from those on the surfaces $\partial N_\rho$
as $\rho\to\infty$ by a simple rescaling, see the main text. The relation between the data
at infinity and those on $\partial N$ are as follows
\be\label{intr-tr}
I^* = \frac{1}{2}(I+\II)I^{-1}(I+\II), \qquad
\II^* = \frac{1}{2}(I+\II)I^{-1}(I-\II).
\ee
These relations can be inverted, with the inverse relations looking exactly the same,
with starred quantities replaced by non-starred everywhere, see (\ref{metric-rel}) below. We note
that the metric $I^*$ is in the conformal class of the boundary $\partial M$, but its
precise form of course depends on the convex domain $N$ used.

\paragraph{The new variation formula.}
The formula (\ref{intr-1}) can then be re-written in terms of the variations
$\delta I^*, \delta \II^*$, 
the result being
\be\label{intr-2}
\delta W(N) = - \frac{1}{4} \int_{\partial N} \langle\delta \II^* - \frac{H^*}{2}
\delta I^*, I^*\rangle  da^*. 
\ee
This formula could be compared to a similar one given by Anderson in dimension
4 \cite{anderson-L2}. 

Thus, the variational formula for the $W$-volume is essentially invariant
under the 
transformation (\ref{intr-tr}). One simple corollary of this formula is that
the extremum of the $W$-volume under variations
of the metric $I^*$ that keep the area of $\partial M$ as defined by $I^*$
fixed occurs 
for the hyperbolic (i.e. constant negative curvature) $I^*$. As there is always the
unique canonical such $I^*$ (of given area), the extremal $W$-volume becomes truly a
functional of the conformal structure of the boundary components of $M$. By considering the
second variation, we also prove that this extremum of $W$ is a maximum.

It is important to note that the second fundamental form $\II^*$ at infinity 
is completely determined by $I^*$ on all the boundary component. This is essentially a consequence
of Bers' simultaneous uniformization. We give a direct proof of this fact in the main text. 
This fact implies that $W(N)$ is a functional of $I^*$ only, a very important
property of the renormalized volume. If one wishes, one can obtain a more canonical functional
that depends on $M$ only, by taking the metric $I^*$ to be the canonical metric
of curvature $-1$ in the conformal class of infinity of $M$. It is this
``extremal'' $W$-volume, called $W_M$ in the latter sections here, 
that is of most interest due to the following.

An important immediate corollary of (\ref{intr-2}) is the theorem by Takhtajan 
and Teo \cite{Takhtajan:2002cc}
that the extremal renormalized volume is equal to the K\"ahler potential on 
the moduli space of Kleinian manifolds. We refer to the main text for a proof of this.
We note that our proof is entirely geometrical and avoids a reasonably complicated cohomology
machinery that is necessary in \cite{Takhtajan:2002cc}. For this reason our
proof can be immediately extended even to situations where the methods of
\cite{Takhtajan:2002cc} are inapplicable, such as manifolds with cone
singularities. See more remarks on this case below.

\paragraph{Positivity.}
The (extremal, i.e the one for the hyperbolic $I^*$) W-volume coincides with a 
multiple of the potential studied by Teo in \cite{Teo}. The result
of this reference implies that the extremal W-volume is a non-negative function on the moduli space
of manifolds, attaining the zero value only on the Fuchsian manifolds. Here we obtain
a similar result for the $W$-volume of the convex core $C_M$ of $M$, see section 4.
These positivity results for the $W$-volume of two different convex domains in
$M$ 
lead us to suggest that the $W$-volume of any convex domain might be positive. 
We do not attempt to answer this question in the present work leaving it for
future research. 

\paragraph{Manifolds with particles and the Teichm\"uller theory of surfaces
  with cone singularities.} 

One key feature of the arguments presented in this work is that they are
always local, 
in the sense that they depend on local quantities defined on the boundaries of
compact subsets of quasi-Fuchsian manifolds. Thus, we make only a very limited
use 
of the fact that the quasi-Fuchsian manifolds are actually quotients of 
hyperbolic 3-space by a group of isometries. One place where this is used is
in the proof of the fact that $\II^*$ is determined by $I^*$ (actually a
direct consequence of the Bers double uniformization theorem).
We expect that all the results should extend from quasi-Fuchsian
(more generally geometrically finite) 
manifolds to the ``quasi-Fuchsian manifolds with particles''  
which were studied e.g. in \cite{minsurf,qfmp}. Those are actually
cone-manifolds, with 
cone singularities along infinite lines running from one connected component of
the boundary at infinity to the other, along which the cone angle is less than
$\pi$. 

One problem towards such an extension is that although in (non-singular)
quasi-Fuchsian setting the Bers double uniformization theorem shows that
everything is determined by the conformal structure at infinity, there is as
yet no such result in the corresponding case ``with particles''. It appears
likely, however, that such a statements holds for ``quasi-Fuchsian manifolds
with particles''; a first step towards it is made in \cite{qfmp}, while the
second step is one of the objects of a work in progress between the second
author and C. Lecuire. 

The result of \cite{qfmp} could actually already be used --- even without a
global Bers type theorem for hyperbolic manifolds with particles --- to obtain
results on the Teichm\"uller-type space of hyperbolic metrics with $n$ cone
singularities of prescribed angles on a closed surface  of genus $g$. Note
that this space, which can be denoted by $\cT_{g,n,\theta}$ (with $\theta=
(\theta_1, \cdots, \theta_n)\in (0,\pi)^n$) is topologically the same as the
``usual'' Teichm\"uller space $\cT_{g,n}$ of hyperbolic metrics with $n$ cusps
(with a one-to-one correspondence from \cite{troyanov}) but it has a natural
``Weil-Petersson'' metric which is different. It might follow from the
considerations made here, extended to quasi-Fuchsian manifolds with particles,
that this ``Weil-Petersson'' metric is still K\"ahler, with the renormalized
volume playing the role of a K\"ahler potential. A global Bers-type theorem
would not be necessary for this because, given any hyperbolic metric $h\in
\cT_{g,n,\theta}$ on a surface $\Sigma$, we can consider the ``Fuchsian''
hyperbolic manifold with particles defined as the warped product
$$ M:= (\Sigma \times \R, dt^2 + \cosh(t)^2 h)~. $$
Clearly the conformal structure at infinity on both connected components of
the boundary at infinity of $M$ are given by $h$. Moreover it follows from
\cite{qfmp} that if $h_-:=h$ and $h_+$ is in a small neighborhood $U\subset
\cT_{g,n,\theta}$ of $h$ 
then there exists a unique quasi-Fuchsian manifold with particles, close to
$M$, with conformal structures at infinity given by $h_-$ and $h_+$. 
The arguments developed here (extended to this singular context) should
show that the renormalized volume is a K\"ahler potential for the natural
Weil-Petersson metric on $\cT_{g,n,\theta}$ restricted to $U$.
We leave such an extension to quasi-Fuchsian cone 
manifolds to future work.

\paragraph{Acknowledgment.}

We would like to thank Rafe Mazzeo for useful comments on a previous version
of this text, and in particular for pointing out reference
\cite{epstein-duke}. 

\section{Preliminaries}

In this section we collect various background materials useful further on in
the paper.

\paragraph{Extrinsic invariants of surfaces in $H^3$.}

Let $S$ be a smooth surface in $H^3$. Its Weingarten (or shape) operator is a
bundle morphism $B:TS\rightarrow TS$ defined by:
$$ Bx := -\nabla_xN~, $$
where $N$ is the unit normal vector field to $S$ and $\nabla$ is the
Levi-Civit\`a connection of $H^3$. $B$ is then self-adjoint with respect to
the induced metric on $S$, which we call $I$ here. The second fundamental form
of $S$ is then defined by:
$$ \II(x,y) := I(Bx,y)=I(x,By)~, $$
for any two vectors $x$ and $y$ tangent to $S$ at the same point. The third
fundamental form of $S$ is defined as: 
$$ \III(x,y):=I(Bx,By)~. $$
When $B$ has no zero eigenvalue, $\III$ is a Riemannian metric on $S$.

\paragraph{The Gauss and Codazzi equations.}

The Weingarten operator satisfies two equations on $S$, the Codazzi equation: 
$$ d^\nabla B=0~, $$
and the Gauss equation: 
$$ \det(B) = K+1~, $$
where $K$ is the curvature of the induced metric $I$ on $S$. $\det(B)$ is
called the {\it extrinsic curvature} of $S$, denoted by $K_e$. 

The Gauss and Codazzi equation are the only relations satisfied by the
first and second fundamental forms of a surface. This can be expressed in a
mildly complicated way as the ``fundamental theorem of surface theory'' stated
below. Here and elsewhere we use a fairly natural convention: given two
bilinear symmetric forms $g$ and $h$ on $S$, with $g$ positive definite, we
denote by $g^{-1}h$ the unique bundle morphism $b:TS\rightarrow TS$,
self-adjoint for $g$, such that $h(x,y)=g(bx,y)$ for any two vectors $x,y$
tangent to $S$ at the same point. For instance, $I^{-1}\II=B$ by definition.

\begin{thm} \label{thm:surf}
Let $g$ and $h$ be two smooth symmetric bilinear forms on a simply
  connected 
  surface $S$, with $g$ positive definite at each point. Define $b:=g^{-1} h$,
  $h={\rm Tr}(b)$, and $k_e={\rm det}(B)$. If $g, h$ satisfy the constraints:
\be \label{Codazzi} \nabla^g h = 0 \qquad {\rm (Codazzi)} \\
\label{Gauss} k_e = K_g+1 \qquad  {\rm (Gauss)},
\ee
where $K_g$ is the Gauss curvature of $g$, 
then there exists a unique immersion of $S$ into the hyperbolic space $H^3$
such that $g, h$ are 
the induced metric and second fundamental forms of $S$ respectively.
\end{thm}

The next lemma, which is elementary and well-known, describes the behavior of
the surfaces at constant distance from a fixed surface.

\begin{lemma} \label{lm:base}
Let $S$ be a surface in $H^3$, with bounded principal curvatures, 
and let $I, B$ be the first fundamental 
form and the shape operator of $S$ correspondingly. Let $S_\rho$ be the
surface at  
distance $\rho$ from $S_\rho$. Then, for sufficiently small $\rho$ the induced
metric on $S_\rho$ is: 
\be \label{metric} I_\rho(x,y) =
I((\cosh(\rho)E+\sinh(\rho)B)x,(\cosh(\rho)E+\sinh(\rho)B)y)~. 
\ee
Here $E$ is the identity operator.  
\end{lemma}

Note that this lemma also holds for a surface $S$ in any hyperbolic 3-manifold
$M$, not necessarily 
$H^3$. We also note that when the surface $S$ is convex, then the expression
(\ref{metric}) gives 
the induced metric on any surface $\rho>0$, where $\rho$ increases in the
convex direction. A proof of this lemma, and of the two corollaries which
follow, can be found in \cite{minsurf}.

\begin{cor} The area of the surfaces $S_\rho$ is given by:
\be\label{area} 
A(\rho) & = & \int_S \det(\cosh(\rho)E+\sinh(\rho)B) da \\
& = & \int_S \left(\cosh^2(\rho) + \cosh(\rho)\sinh(\rho) H + \sinh^2(\rho) K_e \right) da.
\ee
\end{cor}

\begin{cor} The integrated mean curvature of $S_\rho$ is given by:
\be\label{H}
\int_{S_\rho} H da = \frac{\dr }{\dr \rho} A(\rho) = \int_{S} 
(\sinh(2\rho) (1+K_e) + \cosh(2\rho) H) da. 
\ee
\end{cor}

\section{The new volume}

\paragraph{Definition.}

We define the $W$-volume of a compact hyperbolic manifold with boundary as
follows.  

\begin{df} Let $M$ be a hyperbolic 3-manifold, and $N$ be a compact subset of
  $M$ with boundary $\dr N$.  
Then:
\be\label{vol}
W(N) := V(N) - \frac{l^2}{4} \int_{\dr N} H da, 
\ee
where $V(N)$ is the volume of $N$, $l$ is the radius of curvature of the
space,  
$H$ is the mean curvature of the boundary, and $da$ is the
area element of the metric induced on $\dr N$. 
\end{df}

Note that the volume defined
above {\it does not} coincide with the usual Einstein-Hilbert action:
$$ S_{EH}[g]= \int_M \sqrt{g} (R-2\Lambda) + 2\int_{\dr M} H$$
evaluated on a metric of constant curvature. Indeed, for such a metric 
$R-2\Lambda=4\Lambda$. Thus, defining the radius of curvature as
$\Lambda=-1/l^2$ we 
have:
\be\label{EH}
I_{EH}(N):=  -\frac{l^2}{4} S_{EH}(N) = V(N) - \frac{l^2}{2} \int_{\dr N} H,
\ee
which is different from (\ref{vol}). We will set the radius of curvature $l$ to one 
in what follows. 

The following property of $W(N)$ is obvious:
\begin{lemma} 
The $W$-volume is additive: if $N_1$ and $N_2$ are two compact
sub-manifolds of $M$ such that $N_1\cap N_2$ is a disjoint union of connected components
of the boundary of both $N_1$ and $N_2$ then:
$$ W(N_1 \cup N_2) = W(N_1) + W(N_2)~.$$
\begin{proof} This is obvious when $N_1, N_2$ do not share any  boundary
components. For 
$N_1, N_2$ such that a part of their respective boundaries is shared the additivity follows from 
the fact that the mean curvatures of that boundary component have the opposite sign as viewed from
$N_1$ and $N_2$. 
\end{proof}
\end{lemma}
We note that the on-shell Einstein-Hilbert action (\ref{EH}) is also additive.
Thus, at this stage there is no reason to prefer (\ref{vol}) to (\ref{EH}). However,
as we shall see in the next section, it is the quantity (\ref{vol}) that 
behaves much more regularly for non-compact hyperbolic manifolds as well as
for compact ones. Also, as we shall presently see, it is the $W$-volume
that is the self-dual one.

\paragraph{Self-duality.}
Here we prove self-duality of the new $W$-volume by considering its Legendre transform.
We will need to use the variation formula (\ref{intr-1}) given in the introduction.
 An immediate consequence of this formula is that the variation of the $W(N)$ under the condition
that the ``conjugate'' momentum
\be\label{mom}
\pi = - \frac{1}{4} (\II- \frac{H}{2} I)
\ee
is fixed is given by
\be
\delta W(N) = \int_{\partial N} \langle\pi, \delta I\rangle .
\ee
The quantity (\ref{mom}) is exactly the unique combination of $I,\II$ such that when it
is kept fixed the variation of $W(N)$ produces exactly $\pi$. The dual W-volume can now be
obtained by a Legendre transform:
\be
{}^*W(N):= - \int_{\partial N} \langle\pi, I\rangle  + W(N).
\ee
We see that, because $\pi$ is traceless, the $W$-volume is self-dual: ${}^*W(N)=W(N)$.
Note that a similar argument applied to the usual volume $V(N)$ of a domain $N$ shows that
its Legendre transform is given by the Einstein-Hilbert functional $I_{EH}(N)$. This demonstrates
the fact that the Einstein-Hilbert functional is the dual volume $I_{EH}(N)={}^*V(N)$ to which
we referred to in the introduction.

We have so far discussed the duality only in the context of a compact domain $N\subset H^3$.
The duality is, however, more general and holds also for domains in more general hyperbolic
3-manifolds. Note, however, that in this more general context one has to be careful about 
the geometrical meaning of the dual volume. Indeed, the 3-manifold dual to $M$ is modeled
on $dS_3$ spacetime, and typically has two disconnected components, each having an internal
boundary -- a surface dual to the boundary of the convex core in $M$. Given a convex domain
$N$ that contains the convex core $C_M$ one can meaningfully talk about the dual domain
in $dS_3$ as being a domain in the dual manifold located between the surfaces dual to
$\partial N$ and the internal boundary of ${}^*M$. This discussion serves as a good
introduction to the following section in which we discuss precisely those more general
hyperbolic 3-manifolds for which the notion of the $W$-volume is of interest.

\section{Convex co-compact hyperbolic manifolds}

\paragraph{Definitions and first properties.} 

We first need to define convex co-compact hyperbolic manifolds. 

\begin{df} \label{def:man} 
A complete hyperbolic 3-manifold $M$ 
is convex co-compact if there is a compact subset
$N\subset M$ whose boundary $\dr N$ is convex and such that the normal
exponential map from 
$\dr N$ to the conformal boundary $\dr M$ is a homeomorphism. Each connected
component of 
the complement $M\backslash N$ is called a hyperbolic end of $M$.
\end{df}

Note that the condition on $N$ is equivalent to the fact that $N$ is {\it
  strongly convex} 
in the sense that any geodesic segment in $M$ with endpoints in $N$ is
  actually contained in $N$. 

Simple examples of convex co-compact manifolds are: 
Schottky manifolds, each having one hyperbolic end;
quasi-Fuchsian manifolds with two hyperbolic ends.

The new volume $W(N)$ is especially interesting because it is (almost) defined not only for 
compact subsets $N$ as in the above definition, but also for the hyperbolic ends. 
The following computation is central to motivate the definition that follows.

\begin{lemma} Consider a hyperbolic end of a convex co-compact manifold $M$, and let $S$ be a connected
component of the boundary $\dr N$ of the compact subset $N$ of the definition \ref{def:man}. 
The $W$-volume of the sub-manifold contained between the surfaces $S$ and $S_\rho$ is given by:
\be\label{sandwitch}
W[S,S_\rho]= -\frac{\rho}{2} \int_S K da = 2 \pi \rho (g-1),
\ee
where $g$ is the genus of $S$. 
\end{lemma}  

\begin{proof}
We have:
\begin{eqnarray*}
V(\rho) & = & \int_0^\rho A(r)dr \\
& = & \int_0^\rho \int_S \det(\cosh(r)E+\sinh(r)B) da dr \\
& = & \int_0^\rho \int_S \cosh^2(r) + H\cosh(r)\sinh(r) + K_e \sinh^2(r) da dr \\
& = & \int_0^\rho \int_S (\cosh^2(r)+\sinh^2(r)) + K\sinh^2(r) + H\cosh(r)\sinh(r)
da dr \\ 
& = & \int_0^\rho \int_S \cosh(2r) + K\frac{\cosh(2r)-1}{2} + H \frac{\sinh(2r)}{2} da dr~,
\end{eqnarray*}
so that:
\be\label{vol-end}
V(\rho) = \frac{1}{2} \int_S \left( \sinh(2\rho) + \frac{K}{2}(\sinh(2\rho)-2\rho) 
+ \frac{H}{2}(\cosh(2\rho)-1) \right) da.
\ee
Now the $W$-volume is given by:
$$
V(\rho)-\frac{1}{4} \int_{S_\rho} Hda + \frac{1}{4} \int_S Hda.
$$
The formula (\ref{sandwitch}) follows by combining (\ref{vol-end}) with (\ref{H}).
\end{proof}

\paragraph{The relative $W$-volume.}

Thus, the $W$-volume of a portion of a
hyperbolic end is just the thickness of the portion considered times 
$2 \pi(g-1)$. This fact motivates the
following definition:

\begin{df} Let $M$ be a convex co-compact hyperbolic 3-manifold with one or more hyperbolic ends.
Let $N$ be a compact convex subset of $M$ as in the definition \ref{def:man}.
The $W$-volume $W(M,N)$ of $M$ relative to $N$ is defined as
the $W$-volume of N.  
\end{df}

Note that $W(N)$ is the same as the $W$-volume of $M$ with the $W$-volumes (\ref{sandwitch}) of the
hyperbolic ends removed. Thus, one could also refer to the volume $W(M,N)$ as the
{\it renormalized} volume. Indeed, it has a close relation to the renormalized volume that
has appeared in the literature. Here we would like to give a comparison to the renormalized volume.

\begin{df}
Let $M$ be a convex co-compact
hyperbolic 3-manifold, and let $S_r^i$ be a foliation
by equidistant surfaces near each component $i$ of the boundary. 
The renormalized volume of $M$ relative to foliations $S_\rho^{i}$ is defined as 
$$ V_R(M, S_\rho) := \lim_{\rho \rightarrow \infty} V(\rho) - \sum_i \frac{1}{2} A_\rho^i -
\sum_i 2 \pi \rho_i (g_i-1)~, $$
Here $V(\rho)$ is the volume of the subset of $M$ bounded by surfaces $S_\rho^i$, $A_\rho^i$ and $g_i$
are the areas and genera of the surfaces $S_\rho^i$ correspondingly. The limit is a multiple limit
of all $\rho_i\to\infty$.
\end{df}

When $M$ is convex co-compact there is a natural foliation of each end by
surfaces equidistant to the (strongly) convex subset $N$. 
In this case we will talk about the renormalized volume $V_R(M,N)$ 
of $M$ relative to the (strongly) convex subset $N$.

\begin{lemma}
Let $M$ be convex co-compact, and let $N$ be a (strongly) convex subset
of $M$ (as in the definition \ref{def:man}). 
The renormalized volume of $M$ relative to $N$ is the $W$-volume minus a
multiple of the 
Euler characteristic of the boundary:
\be\label{ren-vol}
V_R(M, N) = W(M,N) - \sum_i \pi(g_i-1)~. 
\ee
The sum is taken over the boundary components.
\end{lemma}

This formula is the 3-dimensional case of a formula given by C. Epstein for
the renormalized volume of hyperbolic manifolds in \cite{epstein-duke}. We
include a proof for the reader's convenience.

\begin{proof}
Consider one of the hyperbolic ends. Let $S$ be the corresponding
boundary component of $\dr N$. The area (\ref{area}) of $S_\rho$ can be
rewritten as: 
\be
A(\rho) =  \int_S \left( \cosh(2\rho) + \frac{H}{2} \sinh(2\rho) 
+ \frac{K}{2} (\cosh(2\rho)-1) \right) da.
\ee
Subtracting half of this from (\ref{vol-end}) we get:
$$ V(\rho) - A(\rho)/2 - 2 \pi\rho(g-1) = 
- \frac{1}{4} e^{-2\rho} \int_{S} \left( 2 - H + K \right) da 
- \frac{1}{4} \int_S Hda - \pi(g-1). $$
The result now follows by taking the limit $\rho\to \infty$ and
adding the terms corresponding to all the different ends.
\end{proof}

Thus, the result (\ref{ren-vol}) shows that the renormalized volume relative to a 
strongly convex subset $N$
is basically the volume (\ref{vol}) we have defined, apart from a constant term 
proportional to the Euler characteristic of $\dr N$. However, as it is clear from
the proof, the two quantities agree only after the limit is taken. For a finite
$\rho$ the renormalized volume is a complicated functional, and it is only in 
the limit that a simplification occurs. The volume (\ref{vol}) we have defined 
is in contrast simple even for a finite $\rho$, as well as for any
compact domain $N$. One could try to define an analog of the renormalized
volume for a general subset $N\subset M$ by taking the volume minus half of the
area of $\dr N$. This functional however fails to be additive and is thus of a
very 
limited interest, apart from the limiting case when the surface $\dr N$ is
sent to  
infinity. All this makes it clear that the functional $W(N)$ is much more
natural 
to consider than the one that plays a role in the definition of the
renormalized volume. 

Having motivated and defined the $W$-volume, the natural question to ask
is what this quantity depends on. From its definition one may expect
that it depends on the shape of the convex subset 
$N$ in $M$ in a complicated way.
However, as we shall see, this dependence is actually rather simple in that
the $W$-volume is just a certain functional of the so-called ``asymptotic''
metric constructed using the fundamental forms of the boundary of $N$.
We deal with this in the next section. However, before we study this question,
let us demonstrate certain positivity properties of the $W$-volume.

\paragraph{Positivity estimates on $W$.}
We note that some of the quantities considered
here are always positive on the convex core of a quasi-Fuchsian
hyperbolic manifold. We actually prove this here under a technical hypothesis
which is conjecturally always satisfied.

\begin{lemma}
Let $M$ be a quasi-Fuchsian manifold, and let $C_M$ be its convex
core. Suppose that $C_M$ is the Gromov-Hausdorff limit of a sequence
of convex cores of hyperbolic manifolds with bending laminations along
closed curves. 
Then $I_{EH}(C_M)\geq 0$, with equality exactly when $M$ is
Fuchsian.
\end{lemma}

A short explanation on the hypothesis is needed. Given $M$, let $\lambda$ be
the measured bending lamination on the boundary of its convex core. It is
known (see \cite{bonahon-otal,lecuire}) that $\lambda$ is the limit of the
measured bending laminations of a sequence of quasi-Fuchsian manifolds $M_n$
for which the convex core converge, in the Gromov-Hausdorff distance, to 
the convex core of a quasi-Fuchsian manifold $M'$, and that the measured 
bending lamination on the boundary of the convex core of $M'$ is $\lambda$.
According to a conjecture of Thurston, this last point should imply that
$M'=M$, and then the hypothesis made in the lemma would be useless.

\begin{proof}
We consider in the proof that $M$ is not Fuchsian, since in that
case it is quite obvious that $I_{EH}(C_M)=0$.

Suppose first that the bending lamination of $C_M$ is along
disjoint closed curves. Let $l_i$ and $\lambda_i$ be the lengths
and bending angles at those closed curves. It is then known 
\cite{bonahon-otal} that there exists a one-parameter family
of quasi-Fuchsian manifolds, $M_t, 0\leq t\leq 1$, with $M_0$
Fuchsian, $M_1=M$, and such that, for all $t\in [0,1]$, the
bending lamination of the convex core of $M_t$ is $t$ times
the measured bending lamination of the convex core of $M$.
In other terms, the bending angle of the curve $i$ on the
boundary of the convex core of $M_t$ is $t\lambda_i$. Let
$l_i(t)$ be the length of this curve.

Now a simple computation using the Schl\"afli formula shows that:
$$ \frac{dI_{EH}(C_{M_t})}{dt} = 
- \frac{1}{2} \sum_i t\lambda_i \frac{dl_i(t)}{dt}~. $$
However Choi and Series \cite{choi-series} have recently
proved that, in this context, the matrix of the differential
of the lengths with respects to the angles is negative definite.
Since the $\theta_i(t)=t\lambda_i$ here, it follows that:
$$ \sum_i \lambda_i \frac{dl_i(t)}{dt} <0~,$$
and therefore $I_{EH}(C_{M_t})$ is a strictly
increasing function of $t$. Since it vanishes for $t=0$,
it follows that $I_{EH}(M)>0$.

If the bending lamination of $M$ is general -- i.e., it is not
supported on closed curves -- the result can be obtained, thanks
to the technical hypothesis in the lemma, by
approximating $M$ by quasi-Fuchsian manifolds for which the 
bending lamination of the boundary of the convex core is along
closed curves. 
\end{proof}

\begin{cor}
Under the same hypothesis, $W(C_M)\geq 0$, with equality exactly
when $M$ is Fuchsian.
\end{cor}
\begin{proof}
This immediately follows from the previous lemma since 
$$ W(M)=I_{EH}(M) + \frac{1}{4}\int_{\dr M} H~. $$  
\end{proof}

These results should be compared to a recent result by Teo \cite{Teo} that
shows that the $W$-volume 
extremized in a certain way, to be explained below, is positive. This other
extremized volume 
is simply the volume of a different convex domain in $M$. Taken together, these
results suggest that the $W$-volume of any convex domain $N$ might be positive.
We will not attempt to prove this statement in the present work.

\section{Description ``from infinity''}

\paragraph{The metric at infinity.}

In this section we switch from a description of the renormalized volume from
the boundary of a convex subset to the boundary at infinity of $M$. This
description from infinity is remarkably similar to the previous one from the
boundary of a convex subset.

\begin{lemma}
Let $M$ be a convex co-compact hyperbolic 3-manifold, and let $N\subset M$ be
compact and ``strongly'' convex with smooth boundary. Let $S_\rho$ be the
equidistant surfaces from $\dr N$. The induced metric on $S_\rho$ 
is asymptotic, as $\rho\rightarrow \infty$, to $(1/2) e^{2\rho} I^*$, where
$I^*=(1/2)(I+2\II+\III)$ is defined on $\dr N$.
\begin{proof}
Follows from Lemma \ref{lm:base}.
\end{proof}
\end{lemma}
It is the metric $I^*$ that will play such a central role in what follows, so
we would like to state some of its properties.

\begin{lemma}
The curvature of $I^*$ is 
\be\label{curv*}
K^*:=\frac{2K}{1+H+K_e}~. 
\ee 
\end{lemma}

\begin{proof}
The Levi-Civit\`a connection of $I^*$ is given, in terms of the Levi-Civit\`a
connection $\nabla$ of $I$, by:
$$ \nabla^*_xy = (E+B)^{-1} \nabla_x((E+B)y)~. $$
This follows from checking the 3 points in the definition of the Levi-Civit\`a 
connection of a metric:
\begin{itemize}
\item $\nabla^*$ is a connection.
\item $\nabla^*$ is compatible with $I^*$.
\item it is torsion-free (this follows from the fact that $E+B$ 
verifies the Codazzi equation: $(\nabla_x (E+B))y=(\nabla_y(E+B))x$). 
\end{itemize}

Let $(e_1, e_2)$ be an orthonormal moving frame on $S$ for $I$, and let
$\beta$ be its connection 1-form, i.e.:
$$ \nabla_xe_1 = \beta(x)e_2, ~\nabla_xe_2 = -\beta(x)e_1~. $$
Then the curvature of $I$ is defined as: $d\beta = -K da$. 

Now let $(e^*_1,e^*_2) := \sqrt{2}((E+B)^{-1}e_1,(E+B)^{-1}e_2)$; clearly it
is an orthonormal moving frame for $I^*$. Moreover the expression of
$\nabla^*$ above shows that its connection 1-form is also $\beta$.
It follows that $Kda = -d\beta = K^*da^*$, so that:
$$ K^* = K \frac{da}{da^*} = \frac{K}{(1/2)\det(E+B)} = \frac{2K}{1+H+K_e}~. $$
\end{proof}

We note that the metric $I^*$ is defined for any surface $S\subset
M$. However, it might have singularities (even when the surface $S$ is smooth)
unless $S$ is strictly {\it horospherically convex}, i.e., 
its principal curvatures are less than $1$ (which implies that it 
remains on the concave
side of the tangent horosphere at each point). If $S$ is a strictly 
horospherically convex surface $S$ embedded in
a hyperbolic end of $M$ then the metric $I^*$ 
is guaranteed to be the in the conformal class of the  
(conformal) boundary at infinity of $M$. For
a general surface $S$ the ``asymptotic'' metric has nothing to do with the
conformal infinity, and in particular, does not have to be in the conformal
class of the boundary. 

\paragraph{The $W$-volume as a functional of $I^*$.}
We claim that the $W$-volume of a convex co-compact manifold 
$M$ relative to a convex subset 
$N$ is a functional of only the metric $I^*$ that
is built using the fundamental forms of the boundary $\dr N$. This
claim can be substantiated in several ways. One way is to refer to
the results about the renormalized volume. It is known from
\cite{Krasnov:2000zq}, \cite{Takhtajan:2002cc} 
that the renormalized volume of $M$ is given by
the so-called Liouville functional for the asymptotic metric $I^*$. 
To prove this result one uses an explicit foliation
of the covering space $H^3$ by certain equidistant surfaces. The easy
part of the computation is then to identify the ``bulk'' part of the
Liouville action. The hard part is to show that all the boundary terms
that arise are exactly the boundary of the fundamental domain terms 
necessary to define the Liouville action. As we have said, this
computation is done in the covering space and is not particularly 
illuminating as far as the geometry of the problem is concerned. In
this paper we would like to give a more geometric perspective. We 
demonstrate the above assertion by proving an explicit variational
formula for $W(N)$ in terms of $I^*$. However, before we do this,
we would like to introduce some other quantities defined 
``at infinity''.

\paragraph{Second fundamental form at infinity.}

We have already defined the metric ``at infinity''. Let us now add to this
a definition of what can be called the second fundamental form at infinity.

\begin{df} Given a surface $S$ with the first, second and third
fundamental forms $I,\II$ and $\III$, we 
define the first and second fundamental forms ``at infinity'' as:
\be\label{ff*}
I^* = \frac{1}{2}(I + 2\II + \III)= \frac{1}{2}(I+\II)I^{-1}(I+\II) =
\frac{1}{2}I((E+B)\cdot, (E+B)\cdot)~, \\
\nonumber 
\II^* = \frac{1}{2}(I-\III)=\frac{1}{2}(I+\II)I^{-1}(I-\II) = 
\frac{1}{2}I((E+B)\cdot, (E-B)\cdot)~.
\ee
\end{df}

It is then natural to define:
\be B^*:=(I^*)^{-1}\II^* = (E+B)^{-1}(E-B)~, \label{eq:def-B*} \ee
and 
$$ \III^*:=I^*(B^*\cdot, B^*\cdot) = I((E-B)\cdot, (E-B)\cdot)~. $$
Note that, for a surface which has principal curvatures strictly 
bounded between $-1$ and $1$, $\III^*$ is also a smooth metric and its
conformal class corresponds to that on the other component of the boundary at
infinity. This is a simple consequence of the Lemma \ref{lm:base} and
the fact that when the principal curvatures are strictly bounded
between $-1,1$ the foliation by surfaces equidistant to $S$ 
extends all the way through the manifold $M$. Such manifolds were
called {\it almost-Fuchsian} in our work \cite{minsurf}.

As before, those definitions make sense for any surface, but it is only for a
convex surface (or more generally for a horospherically convex surface)
that the fundamental forms so introduced are guaranteed to have
something to do with the actual conformal infinity of the space.

\paragraph{The Gauss and Codazzi equations at infinity.}

We also define $H^*:=\tr(B^*)$. The Gauss equation for ``usual'' surfaces in
$H^3$ is replaced by a slightly twisted version.

\begin{remark} \label{rk:H*K*}
$H^*=-K^*$: the mean curvature at infinity is equal to minus 
the curvature of $I^*$.
\end{remark}

\begin{proof}
By definition, $H^* = \tr((E+B)^{-1}(E-B))$. An elementary computation (for
instance based on the eigenvalues of $B$) shows that
$$ H^* = \frac{2-2\det(B)}{1+\tr(B)+\det(B)}~. $$
But we have seen (as Equation (\ref{curv*})) that $K^*=2K/(1+H+K_e)$,
the result follows because, by the Gauss equation, $K=-1+\det(B)$. 
\end{proof}

However, the ``usual'' Codazzi equation holds at infinity.

\begin{remark} \label{rk:codazzi*}
$d^{\nabla^*}B^*=0$. 
\end{remark}

\begin{proof}
Let $u,v$ be vector fields on $\dr_\infty M$. Then it follows from the
expression of $\nabla^*$ found above that:
\begin{eqnarray*}
(d^{\nabla^*}B^*)(x,y) & = & \nabla^*_x(B^*y) - \nabla^*_y(B^*x) - B^*[x,y] \\
& = & (E+B)^{-1}\nabla_x((E+B)B^*y) - (E+B)^{-1}\nabla_y((E+B)B^*x) - B^*[x,y]
\\ 
& = & (E+B)^{-1}\nabla_x((E-B)y) - (E+B)^{-1}\nabla_y((E-B)x) -
(E+B)^{-1}(E-B)[x,y] \\
& = & (E+B)^{-1}(d^\nabla (E-B))(x,y) \\
& = & 0~.
\end{eqnarray*}
\end{proof}

\paragraph{Inverse transformations.}
The transformation $I,\II \to I^*, \II^*$ is invertible. The inverse is
given explicitly by:
\begin{lemma} Given $I^*, \II^*$ the fundamental forms $I, \II$ such that
(\ref{ff*}) holds are obtained as:
\be\label{metric-rel}
I = \frac{1}{2} (I^*+\II^*)(I^*)^{-1}(I^*+\II^*) =
\frac{1}{2}I^*((E+B^*)\cdot, (E+B^*)\cdot))~, \\ \nonumber
\II = \frac{1}{2} (I^*+\II^*)(I^*)^{-1}(I^*-\II^*)=
\frac{1}{2}I^*((E+B^*)\cdot, (E-B^*)\cdot))~.
\ee
Moreover,
\be B = (E+B^*)^{-1}(E-B^*)~. \label{eq:B-inverse} \ee
\end{lemma}

Having an expression for the fundamental forms of a surface in terms of the
one at infinity one can re-write the metric of Lemma \ref{lm:base} induced 
on surfaces equidistant to $S$ in terms of $I^*,\II^*$.

\begin{lemma} The metric (\ref{metric}) induced on the surfaces equidistant to
  $S$  
can be re-written in terms of the fundamental forms ``at infinity'' as:
\be \label{metric-infinity}
I_\rho = \frac{1}{2} e^{2\rho} I^* + \II^* + \frac{1}{2} e^{-2\rho} \III^*~.
\ee
\end{lemma}
This lemma shows the significance of $\II^*$ as being the constant term of the metric.
This lemma also shows clearly that when the equidistant foliation extends all
the way through $M$ (i.e. when the principal curvatures on $S$ are in $(-1,1)$),
the conformal structure at the second boundary component of $M$ is that
of $\III^*=\II^* (I^*)^{-1} \II^*$. Thus, in this particular case of
almost-Fuchsian manifolds, the knowledge of $I^*$ on both boundary
components of $M$ is equivalent to the knowledge of $I^*, \II^*$ near either
component. In other words, $\II^*$ is determined by $I^*$. This statement
is more general and works for manifolds other than almost-Fuchsian.

\paragraph{Fundamental Theorem of surface theory ``from infinity''.}

Let us now recall that the Fundamental Theorem of 
surface theory, Theorem \ref{thm:surf}, states that
given $I,\II$ on $S$ there is a unique embedding of $S$ into the
hyperbolic space. Then (\ref{metric}) gives an expression for the
metric on equidistant surfaces to $S$, and thus describes a hyperbolic manifold
$M$ in which $S$ is embedded, in some neighborhood of $S$. It would be 
possible to state a similar result for hyperbolic ends, uniquely 
determined by $I^*$ and $\II^*$ at infinity. But there is also an
analogous theorem, based on a classical result of Bers \cite{Bers},
in which the first (and only the first) form at infinity
is used. This can be compared with arguments used in \cite{horo}.

\begin{thm} 
Given a convex co-compact 3-manifold $M$, and 
a metric $I^*$ (on all the boundary components of $M$) in the conformal class
of the boundary,  
there is a unique foliation of each end of $M$ by convex equidistant surfaces
$S_\rho\subset M$ 
$(1/2)(I_\rho + 2\II_\rho + \III_\rho)=e^{2\rho} I^*$, 
where $I_\rho, \II_\rho, \III_\rho$ are the fundamental forms of $S_\rho$.
\end{thm}

\begin{remark} Note that one does not need to specify $\II^*$. The first
  fundamental 
form $I^*$ (but on all the boundary components) is sufficient.
\end{remark}

\begin{proof}
The surfaces in question can be
given explicitly as an embedding of the universal cover $\tilde{S}$ of $S$
into the hyperbolic space.  
Thus, let $(\xi,y), \xi>0, y\in C$ be the usual upper half-space model
coordinates of $H^3$. Let us write the metric at infinity as 
\be\label{I*-Liouv}
I^* = e^{\phi} |dz|^2,
\ee
where $\phi$ is the Liouville
field covariant under the action of the Kleinian group giving $M$ on $S^2$. 
The surfaces
are given by the following set of maps: $Eps_\rho: S^2 \to H^3, z\to (\xi,y)$
(here $Eps$ stands for Epstein, 
who described these surfaces in \cite{Eps}):
\be 
\xi = \frac{\sqrt{2} e^{-\rho} e^{-\phi/2}}{1+(1/2)e^{-2\rho} e^{-\phi}
  |\phi_z|^2}, \\ \nonumber 
y = z + \phi_{\bar{z}} \frac{e^{-2\rho} e^{-\phi}}{1+(1/2)e^{-2\rho} e^{-\phi}
  |\phi_z|^2}. 
\ee
As is shown by an explicit computation, the metric induced on the 
surfaces $S_\rho$ is given by
(\ref{metric-infinity}) with 
\be
\II^* = \frac{1}{2}(\theta dz^2 + \bar{\theta} d\bar{z}^2) + \phi_{z\bar{z}} dz d\bar{z}, \\ 
\label{theta}
\theta = \phi_{zz} - \frac{1}{2} (\phi_z)^2.
\ee
Thus, we see that $\II^*$ is determined by the conformal factor in (\ref{I*-Liouv}).
\end{proof}

\begin{remark} This theorem itself implies that the renormalized volume only depends on $I^*$. Indeed,
the foliation $S_\rho$ of the ends does depend only on $I^*$, and this foliation can be
used for regularization and subtraction procedure. Then the fact that the $W$-volume is
essentially the renormalized volume implies that $W$-volume is a functional of $I^*$ only.
In the next section we will find a formula for the first variation of this functional.
\end{remark}

\begin{cor} If the principal curvatures at infinity (eigenvalues of $B^*$) are
  positive the map $Eps_\rho$ is a homeomorphism onto its image for any $\rho$.
\end{cor}

\begin{proof}
We first note that the map $Eps_\rho$ is not always a 
homeomorphism, and the surfaces $S_\rho$ are not
necessarily convex, but for sufficiently large $\rho$ 
both things are true. A condition that guarantees that 
$Eps_\rho$ is a homeomorphism for any $\rho$ is stated 
above. This condition can be
obtained from the requirement that the principal 
curvatures of surfaces $S_\rho$ are in $[-1,1]$. Let us
consider the surface $S:=S_{\rho=0}$ the first and 
second fundamental forms of which are given by 
(\ref{metric-rel}) (this immediately follows from (\ref{metric-infinity})). 
The shape operator of
this surface is then given by $B=(E+B^*)^{-1}(E-B^*)$. It is then 
clear that the principal
curvatures of $S$ are given by $k_i = (1-k_i^*)/(1+k_i^*)$, where the
$k_i^*$ are the ``principal curvatures''
(eigenvalues) of $B^*$. The latter are easily shown to be given by
\be
k^*_{1,2} = e^{-\phi}\left(\phi_{z\bar{z}}\pm \sqrt{\theta\bar{\theta}}\right).
\ee
It is now easy to see that the condition $k_{1,2}\in(-1,1)$ 
is equivalent to the condition 
$k^*_{1,2}>0$. This is a necessary and sufficient condition 
for the foliation by surfaces $S_\rho$
to extend throughout $M$. If this condition is satisfied the map $Eps_\rho$ is 
a homeomorphism for any $\rho$.
\end{proof}

Interestingly, this condition makes sense not only in the quasi-Fuchsian situation
but is more general. Thus, for example, it applies to the Schottky manifolds. But for
the Schottky manifolds with their single boundary component the foliation by equidistant surfaces
$S_\rho$ cannot be smooth for arbitrary $\rho$. It is clear that surfaces must develop singularities
for some value of $\rho$. We therefore get a very interesting corollary: 
\begin{cor}
There is no Liouville field $\phi$ on $\C$ invariant under a Schottky group such that $\phi_{z\bar{z}}$ is greater
than $|\phi_{zz}-(1/2)\phi_z^2|$ everywhere on $\C$. 
\end{cor}
\begin{proof}Indeed, if such a Liouville field existed,
we could have used it to construct a smooth equidistant foliation for arbitrary values of $\rho$,
but this is impossible. 
\end{proof}
Similar statement holds for a Kleinian group with more than two components
of the domain of the discontinuity.

\section{The Schl\"afli formula ``from infinity''}

In this section we obtain a formula for the first variation of the
renormalized volume. 
The computation of this section is a bit technical. Readers not interested in
the details 
are advised to skip this section on the first reading. The result of the
computation 
is given by the formula (\ref{dw}) below.

\paragraph{The Schl\"afli formula.}
As we have seen in the previous sections, the renormalized volume of a convex
co-compact hyperbolic 
3-manifold $M$ can be expressed as the W-volume of any convex domain $N\subset
M$. The W-volume is equal to the 
volume of $N$ minus the quarter of the integral of the mean curvature over the
boundary of $N$.  
Let us consider what happens if one changes the metric in $M$.
As was shown in \cite{sem-era}, the following formula for the variation of the volume holds
\be
2 \delta V(N) = \int_{\partial N} \left( \delta H + \frac{1}{2} \langle\delta I,\II\rangle  \right) da.
\ee
Here $H$ is the trace of the shape operator $B=I^{-1} \II$, and the expression
$\langle A,B\rangle $ stands for
$\tr(I^{-1} A I^{-1} B)$. We can use this to get the
following expression for the variation of the W-volume:
\be\nonumber
\delta W(N) = \frac{1}{2} \int_{\partial N} \left( \delta H + \frac{1}{2} \langle \delta I,\II\rangle  \right) da - 
\frac{1}{4} \int_{\partial N} \delta H da - \frac{1}{4} \int_{\partial N} H
\delta(da)~, \ee
so that
\be \label{sch-1}
\delta W(N) = 
\frac{1}{4} \int_{\partial N} \left( \delta H + \langle \delta I,\II- \frac{H}{2} I\rangle  \right) da~.
\ee
To get the last equality we have used the obvious equality
\be
da' = \frac{1}{2} \tr(I^{-1} \delta I) = \frac{1}{2} \langle \delta I,I\rangle  da.
\ee
The formula (\ref{sch-1}) can be further modified using
\be
\delta H = \delta(\tr(I^{-1} \II)) = -\tr( I^{-1} (\delta I) I^{-1} \II) + \tr(I^{-1} \delta\II) = -\langle \delta I,\II\rangle  + \langle I,\delta \II\rangle .
\ee
We get
\be\label{sch-2}
\delta W(N) = \frac{1}{4} \int_{\partial N} \langle \delta \II - \frac{H}{2} \delta I, I\rangle  da.
\ee
It is this formula that will be our starting point for transformations to express the
variation in terms of the data at infinity.

\paragraph{Parameterization by the data at infinity.}

Let us now recall that, given the data $I,\II$ on the boundary of $N$ one can
introduce the 
first and second fundamental forms ``at infinity'' via
(\ref{ff*}). Conversely, knowing 
the fundamental forms $I^*,\II^*$ ``at infinity'' one can recover the
fundamental forms 
on $\partial N$ via (\ref{metric-rel}). Our aim is to rewrite the variation
(\ref{sch-2}) 
of the W-volume in terms of the variations of the forms $I^*,\II^*$. 

\begin{lemma} \label{lm:schafli-inf}
The first-order variation of $W$ can be expressed as
\be\label{dw}
\delta W(N) = - \frac{1}{4} \int_{\partial N} \langle  \delta \II^* -
\frac{H^*}{2} \delta I^*, I^*\rangle  da^*. 
\ee
\end{lemma}

\begin{proof}
Clearly the first-order variation of $W$ contains two kinds of terms, those
related to the first-order variation of $I^*$ and those coming from the
first-order variation of $B^*$, which we call $\delta B^*$ here. 
We define $X:=(I^*)^{-1}\delta I^*$, and consider separately the terms which
are linear in $X$ and those which are linear in $\delta B^*$. To simplify
notations we use the notation ``$O(X)$'' (resp. ``$O(\delta B^*)$'') to
describe any term linear in $X$ (resp. in $\delta B^*$). We will be using 
repeatedly a simple formula valid for any two $2\times 2$ 
matrices $A$ and $B$:
\be\label{sch-4}
\det(A) \, \tr(A^{-1} B) = \tr(A)\, \tr(B) - \tr(AB)~.
\ee

We consider first the terms linear in $\delta B^*$. We already know that 
$$ 2I = I^*((E+B^*)\cdot, (E+B^*)\cdot)~, $$
it follows that
$$ 2\delta I = 2I^*((E+B^*)\cdot, \delta B^*) + O(X)~, $$
so that
$$ \langle \delta I,I\rangle = 2\tr((E+B^*)^{-1}\delta B^*) + O(X)~. $$
Similarly, we know that
$$ 2\II = I^*((E+B^*)\cdot, (E-B^*)\cdot)~, $$
it follows that
$$ 2\delta \II = I^*(\delta B^*\cdot, (E-B^*)\cdot) - I^*((E+B^*)\cdot, \delta
B^*\cdot) + O(X)~, $$
and therefore
that
$$ \langle \delta \II,I\rangle = \tr((E+B^*)^{-1}(E-B^*)(E+B^*)^{-1}\delta
B^*) - \tr((E+B^*)^{-1}\delta B^*) + O(X)~. $$
Moreover
$$ H = \tr((E+B^*)^{-1}(E-B^*)), $$
and putting all terms together shows that
$$ \left\langle \delta \II - \frac{H}{2}\delta I, I\right\rangle =
- \tr((E+B^*)^{-1}\delta B^*) + $$
$$ + \tr((E+B^*)^{-1}(E-B^*)(E+B^*)^{-1}\delta
B^*) - \tr((E+B^*)^{-1}(E-B^*)) \tr((E+B^*)^{-1}\delta B^*) + O(X)~. $$
The last two terms can be treated as the right-hand side of (\ref{sch-4}),
the equation becomes: 
$$ \left\langle \delta \II - \frac{H}{2}\delta I, I\right\rangle =
- \tr((E+B^*)^{-1}\delta B^*) - \det((E+B^*)^{-1}(E-B^*))
\tr((E-B^*)^{-1}\delta B^*) + O(X)~. $$
We now apply (\ref{sch-4}) to each of the two terms on the right-hand side and
get for the above quantity:
\begin{eqnarray*}
  \frac{\tr((E+B^*)\delta B^*) - \tr(E+B^*) \tr(\delta B^*)
  - \tr(E-B^*)\tr(\delta B^*) + \tr((E-B^*)\delta B^*)}{\det(E+B^*)} + O(X)
\\ 
=  \frac{2\tr(\delta B^*) - 4\tr(\delta B^*)}{\det(E+B^*)} + O(X) =
\frac{-2\tr(\delta B^*)}{\det(E+B^*)}+ O(X) ~. 
\end{eqnarray*}
 
The part which is linear in $X$ can be computed in a similar way. First note
that, by definition of $X$,
$$ \langle \delta I,I\rangle = \tr(X) + O(\delta B^*)~, $$
while 
$$ \delta \II = I^*((E+B^*)\cdot, X(E-B^*)\cdot)+ O(\delta B^*)~, $$
so that 
$$ \langle \delta \II, I\rangle = \tr((E+B^*)^{-1}X(E-B^*))+ O(\delta B^*)~. $$
It follows, using (\ref{sch-4}) and forgetting all terms linear in $\delta
B^*$, that
\begin{eqnarray*}
\langle \delta\II -(H/2)\delta I,I\rangle & = & \tr((E+B^*)^{-1}X(E-B^*))
- (1/2)\tr((E+B^*)^{-1}(E-B^*)) \tr(X) \\
& = & \frac{1}{\det(E+B^*)}(\tr(E+B^*)\tr(X(E-B^*)) - 
\tr((E+B^*)X(E-B^*)) \\
& - & (1/2)\tr(E+B^*)\tr(E-B^*)\tr(X) + (1/2) \tr(E-(B^*)^2) \tr(X)) \\
& = & \frac{1}{\det(E+B^*)}(\tr(E+B^*)\tr(X(E-B^*)) -
\tr((E-(B^*)^2)X) \\
& - & (1/2)(\tr(E)^2-\tr(B^*)^2)\tr(X) + (1/2) \tr(E-(B^*)^2) \tr(X))~.
\end{eqnarray*}
The terms involving $(B^*)^2$ can be replaced using the fact that
$(B^*)^2 - \tr(B^*)B^* + \det(B^*)E=0$, so that $\tr((B^*)^2)=\tr(B^*)^2
-2\det(B^*)$. It follows that, still forgetting all terms which are linear in
$\delta B^*$, we have
\begin{eqnarray*}
\langle \delta\II -(H/2)\delta I,I\rangle & = &
\frac{1}{\det(E+B^*)}(\tr(E+B^*)\tr(X(E-B^*)) -
\tr((1+\det(B^*))X - \tr(B^*)B^*X) \\
& - & (1/2)(4-\tr(B^*)^2)\tr(X) 
+ (1/2)(2-\tr(B^*)^2 +2\det(B^*)) \tr(X)) \\
& = & \frac{1}{\det(E+B^*)}((2+\tr(B^*))(\tr(X)-\tr(B^*X))
- (1+\det(B^*))\tr(X) + \tr(B^*)\tr(B^*X) \\
& + & (\det(B^*)-1)\tr(X)) \\
& = & \frac{\tr(B^*)\tr(X)-2\tr(B^*X)}{\det(E+B^*)}~. 
\end{eqnarray*}

Putting together the terms in $X$ and the terms in $\delta B^*$, and using the
fact that $da^*=(1/2)\det(E+B^*) da$ by definition of $I^*$, we find that:
$$ \langle \delta \II - (H/2)\delta I,I\rangle = -(\tr(\delta B^*) + \tr(B^*X)
- (1/2)\tr(B^*)\tr(X)) (da^*/da) = - \langle \delta \II^* - (H^*/2)\delta I^*,
I^* \rangle (da^*/da)~, $$
and the result clearly follows.
\end{proof}

Formula (\ref{dw}) looks very much like the original formula (\ref{sch-2}),
except for the minus sign 
and the fact that the quantities at infinity are used. The fact that we have
got the same variational 
formula as in terms of the data on $\partial N$ is not too surprising. Indeed,
the variational 
formula (\ref{dw}) was obtained from (\ref{sch-2}) by applying the
transformation (\ref{metric-rel}). As it is clear from  
(\ref{ff*}), this transformation applied twice gives the identity map. In view
of this, it is hard to think of any other 
possibility for the variational formula in terms of $\delta I^*, \delta \II^*$
except being given by the same 
expression (\ref{sch-2}), apart from maybe with a different sign. This is
exactly what we see in (\ref{dw}).  

There is another expression of the first-order variation of $W$, dual to
(\ref{sch-1}), which will be useful below.

\begin{cor} \label{cr:variation}
The first-order variation of $W$ can also be expressed as
$$ \delta W = -\frac{1}{4}\int_{\dr N} \delta H^* + \langle \delta I^*,
\II^*_0\rangle da^*~, $$
where $\II^*_0$ is the traceless part (for $I^*$) of $\II^*$.
\end{cor}

\section{Conformal variations of the metric at infinity}

In this section we use Corollary \ref{cr:variation}
to prove one simple corollary.
Thus, we show that, when varying the W-volume with the area of the boundary
defined by the $I^*$ metric 
kept fixed, the variational principle implies the metric $I^*$ to have
constant negative curvature. 
The variations we consider in this section do not change the conformal
structure of the metric $I^*$, 
and thus do not change the manifold $M$. Geometrically they correspond to
small movements of the 
surface $\partial N$ inside the fixed manifold $M$. 

\paragraph{First variation.} 

We consider in this section a conformal deformation of the metric $I^*$, i.e.,
$\delta I^*=2u I^*$, where $u$ is some function on $\dr N$. Clearly for 
such variations $\langle \delta I^*,\II^*_0\rangle=0$, precisely because
$\II^*_0$ is traceless.

Let us consider the following functional
\be\label{var-3}
F(N) = W(N) - \frac{\lambda}{4} \int_{\partial N} da^* 
\ee
appropriate for finding an extremum of the W-volume with the area computed
using the 
metric $I^*$ kept fixed. The first variation of this functional gives,
using Corollary \ref{cr:variation}:
$$ \delta F = -\frac{1}{4}\int_{\dr N} (\delta H^*) da^* - \frac{\lambda}{4}
\int_{\dr N} 4u da^* = \frac{1}{4}\int_{\dr N} (\delta K^*) da^* -
\frac{\lambda}{4} \int_{\dr N} 4u da^*~. $$
But 
$$ \delta\int_{\dr N}K^* da^*= \int_{\dr N} (\delta K^*) + 4u K^* da^* =0 $$
by the Gauss-Bonnet formula, so that
$$ \delta F = \int_{\dr N} (- uK^* - u\lambda) da^*~. $$
It follows that critical points of $F$ are characterized by the fact that
$K^*=-\lambda$. 

\paragraph{Second variation.}

In this paragraph we would like to verify whether the extremum found above is
maximum or minimum. It is easy to compute the formula for the second
variation,
but for the sake of simplicity we only do it here for conformal deformations
of $I^*$, corresponding to movements of $\partial N$. 

According to the formula found right above, the gradient of $F$ on the 
space of metrics conformal to $I^*$ is simply $DF=-(K^*+\lambda)/4$. 
However a well-known formula on the conformal deformations of metrics
(see e.g. chapter 1 of \cite{Be}) indicates that, under a conformal
deformation $\delta I^*=2vI^*$,
$$ \delta K^* = -2vK^* + \Delta v~, $$
where $\Delta$ is the Laplacian for $I^*$.
It follows directly that the Hessian of $F$ at a critical point is given by
$$ (\mbox{Hess}(F))(2uI^*, 2vI^*)=\int_{\dr N} 2K^*uv-(\Delta v)u da^* =
\int_{\dr N} -2\lambda uv - \langle du,dv\rangle da^*~. $$
This quantity is negative definite (because $\lambda=-K^*>0$ by Gauss-Bonnet)
so that the critical points of $F$ are local maxima.

\section{The renormalized volume as a function on Teichm\"uller space}

In this section we consider the renormalized volume as a function over the
Teichm\"uller space of $\dr N$; in other terms, for each conformal class on
$\dr N$, we consider the extremum of $W$ over metrics of given area within this
conformal class. We have seen in the previous section that this extremum is
obtained at the (unique) constant curvature metric. The main goal here
is to recover by simple differential geometric methods important results of
McMullen \cite{McMullen} -- concerning 
his ``quasi-Fuchsian reciprocity'' -- and 
Takhtajan and Zograf \cite{TZ-Schottky}, 
Takhtajan and Teo \cite{Takhtajan:2002cc} -- 
showing that the renormalized volume provides a K\"ahler
potential for the Weil-Petersson metric. So the ``volume'' that we consider
here is now defined as follows.

\begin{df}
Let $g$ be a convex co-compact hyperbolic metric on $M$, and let $c\in
\cT_{\dr M}$ be the conformal structure induced on $\dr_\infty M$. We call
$W_M(c)$ the value of $W$ on the equidistant foliation of $M$ near infinity
for which $I^*$ has constant curvature $-1$.
\end{df}

In other terms, by the results obtained in the previous section, $W_M(c)$ is
the maximum of $W$ over the metrics at infinity which have the same area as a hyperbolic metric, for
each boundary component of $M$. Throughout this section
the metric at infinity $I^*$ that we consider is the hyperbolic metric, while
the second fundamental form at infinity, $\II^*$, is uniquely determined by
the choice that $I^*$ is hyperbolic. Its traceless part is denoted by
$\II^*_0$. 

\paragraph{The second fundamental form at infinity as the real part of a
  HQD.}

It is interesting to remark that, in the context considered here -- when $I^*$
has constant curvature -- the second fundamental form at infinity has a
complex interpretation. This can be compared with the same phenomenon,
discovered by Hopf \cite{hopf}, for the
second fundamental form of constant mean curvature surfaces in 3-dimensional
constant curvature spaces.

\begin{lemma}
When $K^*$ is constant, $\II^*_0$ is the real part of a quadratic holomorphic
differential (for the complex structure associated to $I^*$) on $\dr_\infty
M$.  This holomorphic quadratic differential is given explicitly by
(\ref{theta}).
\end{lemma}

\begin{proof}
By definition $\II^*_0$ is traceless, which means that it is at each point the
real part of a quadratic differential: $\II^*_0=Re(h)$. Moreover, we have seen 
in Remark \ref{rk:codazzi*} that $B^*$ satisfies the Codazzi equation,
$d^{\nabla^*} B^*=0$. It follows as for constant mean curvature surfaces (see
e.g. \cite{minsurf}) that $h$ is holomorphic relative to the complex structure
of $I^*$.
\end{proof}

\paragraph{The second fundamental form as a Schwarzian derivative.}

The next step is that, for each boundary component $\dr_iM$ of $M$, 
$\II^*_{0i}$ is actually the real part of the Schwarzian
derivative of a natural equivariant 
map between the hyperbolic plane (with its canonical
complex projective structure) to $\dr_iM$ with its complex projective structure
induced by the hyperbolic metric on $M$. In the terminology used by 
McMullen \cite{McMullen}, $\II^*_{0i}$ is the difference between the complex
projective structure at infinity on $\dr_iM$ and the Fuchsian projective
structure on $\dr_iM$.

A simple way to prove this assertion
is to use the formula (\ref{theta}) for the holomorphic quadratic differential
$\theta$ whose real part gives the traceless part of $\II^*$. The Liouville
field $\phi$ that enters into this formula can be simply expressed in
terms of the conformal map from $\dr_iM$ to the hyperbolic plane. It is
then a standard and simple computation to verify that $\theta$ is equal to
the Schwarzian derivative of this map, see e.g. \cite{TZ-Schottky}. To make
this paper self-contained we decided to include yet a different, more geometric proof,
which is spelled out in the appendix. The proof we give is
elementary, based on the conformal factor between the $I^*$
metrics on corresponding surfaces in two foliations. 
It can be compared to the argument used in \cite{McMullen}. 

To state the result, let us 
call $\sigma_F$ the ``Fuchsian'' complex projective structure on $\dr_iM$,
obtained by applying the Poincar\'e uniformization theorem to the conformal
metric at infinity on $\dr_iM$. The universal cover of $\dr_iM$,
with the complex projective structure lifted from $\sigma_F$, is projectively
equivalent to a disk in $\C P^1$. We also call $\sigma_{QF}$ the projective
structure induced on $\dr_iM$ by the hyperbolic metric on $M$. Here ``$QF$''
stands for quasi-Fuchsian (while $M$ is only supposed to be convex
co-compact), this notation is used to keep close to the notation in
\cite{McMullen}. The map $\phi:(\dr_iM,\sigma_F) \rightarrow
(\dr_iM,\sigma_{QF})$ which is isometric for the hyperbolic metrics on
$\dr_iM$, is conformal but not projective between $(\dr_iM,\sigma_F)$
and $(\dr_iM,\sigma_{QF})$, so we can consider its Schwarzian derivative 
$\cS(\phi)$. 

\begin{lemma} \label{lm:schwarzian}
$\II^*_0=-Re(\cS(\phi))$. 
\end{lemma}

It is possible to reformulate this statement slightly by setting
$\theta_i:=\cS(\phi)$ (this is analogous to the notations used in \cite{McMullen}, the
index $i$ is useful to recall that this quantity is related to $\dr_iM$). 
Then $\theta_i$ is a quadratic holomorphic differential (QHD) on $\dr_iM$,
and, still using the notations in \cite{McMullen}, the definition of
$\theta_i$ can be rephrased as: $\theta_i=\sigma_{QF}-\sigma_F$. 
The Lemma can then be written as: $\II^*_{0i}=Re(\theta_i)$. 
A geometric proof of this lemma is given in the appendix.

\begin{remark} {\rm 
Note that $\theta_i$ can also be considered as a complex-valued 1-form on the
Teichm\"uller space of $\dr_iM$. Indeed, it is well known that the cotangent
vectors to $\cT_S$, where $S$ is a Riemann surface, can be described as 
holomorphic quadratic differentials $q$ on $S$. The pairing with a tangent vector (Beltrami differential $\mu$) is given
by the integral of $q\mu$ over $S$. The complex structure on $\cT_S$ can then be described
as follows: the image of the cotangent vector $q$ under the action of the complex structure $J$ is
simply $J(q)=iq$. Another, more geometric way to state the action of $J$ is to note that
it exchanges the horizontal and vertical foliations of $q$. 
Thus, holomorphic quadratic differentials $q$ on $S$ are actually
holomorphic 1-forms on $\cT_S$.}
\end{remark}

\paragraph{The second fundamental form as the differential of $W_M$.}

There is another simple interpretation of the traceless part of the
second fundamental form at infinity.

\begin{lemma} \label{lm:gradient}
The differential $dW_M$ of the renormalized volume $W_M$, 
as a 1-form over the Teichm\"uller space of $\dr M$, 
is equal to $(-1/4)\II^*_0$.
\end{lemma}

\begin{proof}
This is another direct consequence of Corollary \ref{cr:variation} 
because, as one
varies $I^*$ among hyperbolic metrics, $H^*$ (which is equal to $K^*$) remains
equal to $-1$, so that $\delta H^*=0$.  
\end{proof}

\begin{cor}
$\theta_i=-4\partial W_M$.
\end{cor}

\begin{proof}
This follows directly from the lemma, since we already know that $\theta_i$ is
a holomorphic differential.  
\end{proof}

\begin{remark} {\rm 
We would like to emphasize how much simpler is the proof given above
than that given in \cite{TZ-Schottky}, \cite{Takhtajan:2002cc}. Unlike in these
references, which obtain the above result on the gradient of $W_M$ using an
involved computation, the Corollary \ref{cr:variation} implies this result in
one line. This demonstrates the strength of the geometric method used here.}
\end{remark}

\paragraph{McMullen's quasi-Fuchsian reciprocity.}

We can now recover McMullen's quasi-Fuchsian reciprocity as a simple
consequence of the proved above relation between $W_M$ and the 
second fundamental form at
infinity. In the context considered here, it is just a consequence of the fact
that the Hessian of a function is symmetric. 

From this point and until the end of this section we suppose that $M$ is
quasi-Fuchsian. It has two boundary components, which we call $\dr_-M$ and
$\dr_+M$, which are homeomorphic. 
The renormalized volume $W_M$ is now a function on the space of
quasi-Fuchsian metrics on $M$, which by the Bers theorem is
$\cT_{\dr_-M}\times \cT_{\dr_+M}$.

Let $(c_-, c_+)\in \cT_{\dr_-M}\times \cT_{\dr_+M}$, they define a unique
hyperbolic metric on $M$, and therefore a complex projective structure
$\sigma_{QF,+}(c_-,c_+)$ (resp. $\sigma_{QF,-}(c_-,c_+)$) on $\dr_+M$
(resp. $\dr_-M$). Using the Schwarzian derivative construction
this can be used to define a QHD $\beta_+(c_-, c_+)$ (resp. $\beta_+(c_-,
c_+)$) on $\dr_+M$ (resp. $\dr_-M$) as
$$ \beta_+(c_-,c_+) = \sigma_{QF,+}(c_-,c_+)-\sigma_F(c_+)~, $$
and respectively
$$ \beta_-(c_-,c_+) = \sigma_{QF,-}(c_-,c_+)-\sigma_F(c_-)~. $$

Then $\beta_+(c_-,c_+)$ is a HQD on $\dr_+M$, and can therefore be identified
with an element of $T^*_{c_+}\cT_{\dr_+M}$. So $c_-$ determines a cotangent
vector field $\beta_+(c_-, \cdot)$ on $\cT_{\dr_+M}$, and similarly $c_+$
determines a cotangent
vector field $\beta_-(\cdot, c_+)$ on $\cT_{\dr_+M}$.
By Lemma \ref{lm:schwarzian} above, 
$$ \II^*_{0,+} = - Re(\beta_+(c_-,c_+))~, \II^*_{0,-} = 
- Re(\beta_-(c_-,c_+))~. $$

Yet another way to state this relationship is that
\be 
\forall v_+\in T_{c_+}\cT_{\dr_+M}, \beta_+(c_-,c_+)(v_+)=- \langle 
\II^*_{0,+}, v_+ \rangle~, \label{eq:dW} 
\ee
and similarly for $\beta_-$. Using Lemma \ref{lm:gradient}, 
$\beta_+$ and $\beta_-$
can be expressed in terms of $W_M$ as follows:
\be \forall v_+\in T_{c_+}\cT_{\dr_+M}, \beta_+(c_-,c_+)(v_+)=
4dW_M((0,v_+))~,
\label{eq:dW2} \ee
and similarly for $\beta_-$. Here $W_M$ is considered as a function on
$\cT_{\dr_-M}\times \cT_{\dr_+M}$, and $(0, v_+)$ is a vector tangent to
$\cT_{\dr_-M}\times \cT_{\dr_+M}$.

We can now consider the differential of the $\beta_+(\cdot, c_+)$ considered
as a function of the conformal structure on $\dr_-M$, 
$$ D\beta_+(\cdot, c_+):T_{c_-}\cT_{\dr_-M}\rightarrow T^*_{c_+}\cT_{\dr_+M}~, $$
and of $\beta_-(c_-,
\cdot)$, considered as a function of the conformal structure on $\dr_+M$
$$ D\beta_-(c_-, \cdot):T_{c_+}\cT_{\dr_+M}\rightarrow T^*_{c_-}\cT_{\dr_-M}~. $$

\begin{thm}[McMullen's quasi-Fuchsian reciprocity \cite{McMullen}] 
\label{tm:qfr}
The maps $D\beta_+(\cdot, c_+)$ and $D\beta_-(c_-, \cdot)$ 
are adjoint to each other.
\end{thm}

\begin{proof}
Let $v_-\in T_{c_-}\cT_{\dr_-M}$ and $v_+\in T_{c_+}\cT_{\dr_+M}$. Then
$$ \langle D\beta_+(c_-,c_+)(v_-,0), v_+\rangle = 4(D_{(v_-,0)}dW_M)((0,v_+)) = 
4(\mbox{Hess}(W_M))((v_-,0),(0,v_+))~. $$
(Note that the Hessian here can be considered without reference to
a specific connection, because $(v_-,0)$ and $(0,v_+)$ are tangent
to $\cT_{\dr_-M}$ and $\cT_{\dr_+M}$ respectively.)
The symmetry in the right-hand side shows quite clearly that
$$ \langle D\beta_+(c_-,c_+)(v_-,0), v_+\rangle = \langle 
D\beta_-(c_-,c_+)(0, v_+), v_-\rangle~, $$
which is precisely the statement of the theorem.
\end{proof}

\paragraph{The renormalized volume as a K\"ahler potential.}

Finally we show here how to recover in this manner the result of Takhtajan and
Teo \cite{Takhtajan:2002cc} stating that the renormalized volume $W_M$ with
$c_-$ fixed is a K\"ahler potential for the Weil-Petersson metric on
$\cT_{\dr_+M}$. To simplify notations a little, we set
$\theta_{c_-}:=\beta_+(c_-, \cdot)$. 
Since we already know that $\theta_{c_-}=-4\partial W_M$, we
only need to prove that $\db(i\theta_{c_-}) = -2\omega_{WP}$, where
$\omega_{WP}$ is the K\"ahler form of the Weil-Petersson metric on
$\cT_{\dr_+M}$. 

An important part of the argument is that $d\theta_{c_-}$, as a 2-form on
$\cT_{\dr_+M}$, does not depend on $c_-$. This appears as Theorem 7.2 
in McMullen's paper \cite{McMullen}. We include a proof for completeness,
following the proof given in \cite{McMullen}.

\begin{prop}
The differential $d\theta_{c_-}$, considered as a complex-valued 2-form
on $\cT_{\dr_+M}$, does not depend on $c_-$.
\end{prop}

\begin{proof}
Let $v_-\in T_{c_-}\cT_{\dr_-M}$, we want to show that 
the corresponding first-order variation $D_{v_-}(d\theta_{c_-})$
of $d\theta_{c_-}$ vanishes. This will follow from the fact that 
the first-order variation of $\theta_{c_-}$ corresponding to $v_-$,
$D_{v_-}\theta_{c_-}$, is the differential of a function defined 
on $\cT_{\dr_+M}$, namely the function $f_{v_-}$ defined by
$$ f_{v_-}(c_+) = \langle \beta_-(c_-, c_+), v_-\rangle~, $$
where $\langle , \rangle$ is the WP pairing.

The fact that $D_{v_-}\theta_{c_-}=df_{v_-}$ can be proved by
evaluating both sides on a vector $v_+\in T_{c_+}\cT_{\dr_+M}$ and using
the quasi-Fuchsian reciprocity.
$$ \langle D_{v_-}\theta_{c_-}, v_+\rangle 
= \langle D\beta_+(c_-, c_+)(v_-,0), v_+,\rangle 
= \langle D\beta_-(c_-,c_+)(0,v_+), v_-\rangle = 
df_{v_-}(v_+)~. $$
It clearly follows that $d\theta_{c_-}$, as a 2-form on $\cT_{\dr_+M}$,
does not depend on $c_-$.
\end{proof}

That $W_M$ is a K\"ahler potential 
is then reduced to a simple computation in the Fuchsian
situation. 

\begin{prop} \label{pr:varII}
Suppose that $M$ is a Fuchsian manifold, with $c_+=c_-$. Let $I^*$ be the
hyperbolic metric in the conformal class $c_+$. Under a first-order
deformation which does not change $c_-$, 
the variation of $I^*$ and $\II^*_0$ on $\dr_+M$ are related by:
$$ \delta \II^*_0 = - \delta I^*~. $$
\end{prop}

\begin{proof}
It follows from the constructions in the appendix that $B^*=(1/2)E$. Under a
first-order variation we have that
$$ \delta \II^* = (1/2) \delta I^* + I^*(\delta B^*\cdot, \cdot)~, $$
$$ \delta \III^* =  \delta I^*(B^*\cdot, B^*\cdot) + 2\II^*(\delta B^*\cdot,
B^* \cdot) = (1/4)\delta I^* + I^*(\delta B^*\cdot, \cdot)~, $$
so that $\delta\III^* = 2\delta^*\II^*-\delta I^*$. 

But $\III^*$ is the hyperbolic metric in the conformal class $c_-$ and
therefore it does not change under the deformation, so that $2\delta
\II^*=\delta I^*$. 

Moreover, since $I^*$ remains hyperbolic, $H^*=\tr_{I^*}\II^*=1$, 
so that $\II^*_0 = \II^*-(1/2)I^*$, and so 
$$ \delta \II^*_0 = \delta\II^*-(1/2)\delta I^*= - \delta I^*~. $$
\end{proof}

We can reformulate this statement by calling $\theta_R:=Re(\theta_{c_-})$, so
that, by Lemma \ref{lm:schwarzian}, 
$\theta_R(X)=- \langle X,\II^*_0\rangle$. Using the previous proposition,
this can the be stated as
$$ (D_X\theta_R)(Y) = \langle X,Y\rangle_{WP}~, $$
where $D$ is the Levi-Civit\`a connection of the Weil-Petersson metric on
$\cT_{\dr_+M}$. 

We can now compute explicitly an expression of $\db\theta_{c_-}$, denoting by
$J$ the complex structure on $\cT_{\dr_-M}$.
\be 
\db\theta_{c_-}(X,Y) & = & (D_X\theta_{c_-})(Y) + i(D_{JX}\theta_{c_-})(Y) \\
& = & (D_X\theta_R)(Y) - i(D_X\theta_R)(JY) + i((D_{JX}\theta_R)(Y) -
i(D_{JX}\theta_R)(JY)) \\
& = & \langle X,Y\rangle - i \langle X,JY\rangle + i\langle JX,Y\rangle +
\langle JX,JY\rangle \\
& = & 2(\langle X,Y\rangle -i\langle X,JY\rangle)~. 
\ee 

This means precisely that $\db \theta_{c_-}(X,JX)=2i\|X\|^2_{WP}$, and we
recover the result of Takhtajan and Teo \cite{Takhtajan:2002cc} that $W_M$ is
a K\"ahler potential for the Weil-Petersson metric.

\appendix

\section{Appendix: A geometric proof of $\II^*_0=-Re(\cS(\phi))$.}

The proof is based on some preliminary but simple arguments on the conformal
factor relating corresponding surfaces in the equidistant 
foliations of a quasi-Fuchsian and a Fuchsian manifold.

We consider a point $x_0\in \dr_iM$. For each point $x$ in some neighborhood of
$x_0$ there is then a unique geodesic $\gamma_x$,
orthogonal to the leaves of the foliation by the surfaces $\Sigma_\rho$
already introduced above, with endpoint $x$. In addition to $M$, we also
consider the ``Fuchsian'' manifold $M'$ with boundary at infinity equal to the
disjoint union of two copies of $\dr_iM$ (with its conformal metric). So $M'$
is topologically $\dr_iM\times \R$, with the hyperbolic metric
\be d\rho^2 + \cosh^2(\rho) I^*~, \label{eq:metric} \ee
where $I^*$ is the hyperbolic metric on $\dr_iM$. Using this construction we
can 
identify $\dr_iM$ also with the ``upper'' boundary component of $M'$. To avoid
ambiguities we call this upper boundary component $\dr_iM'$, and consider the
map 
$\phi$ defined above as a map between $\dr_iM'$ and $\dr_iM$, and call $x'_0$
the point corresponding to $x_0$ on $\dr_iM'$. There is a
canonical foliation on $M'$, given by the level surfaces $\Sigma'_\rho$
of $\rho$ in (\ref{eq:metric}), to which we can associate a metric $I'^*$ and
a second fundamental form $\II'^*$ at infinity on $\dr_iM'$. Taking the right
choice of foliation near each boundary component (which involves an affine
transformation on $\rho$) leads to $2I'^*=\II'^*=2\phi^*I^*$, so that
$B'^*=(1/2)E$. 

We now consider the geodesics $\gamma_{x_0}$ (in $M$) and $\gamma'_{x'_0}$ (in
$M'$), with their respective parameterizations.
There is a unique hyperbolic isometry $\Phi_0$ sending
$\gamma'_{x'_0}$ to $\gamma_{x_0}$, preserving the parameterization by $\rho$,
and such that the differential at $x'_0$ of the boundary map $\phi_0:\dr_\infty
H^3\rightarrow \dr_\infty H^3$ is tangent to
the isometry between $I'^*$ and $I^*$. Then
$\phi_0$ is a complex projective map but it is {\it not} an isometry between
$I'^*$ and $I^*$ (although its differential at $x'_0$ is isometric). 

We will consider the map $\phi_0^{-1}\circ \phi:\dr_iM'\rightarrow \dr_iM'$,
and show 
that it is tangent to the identity at order 2 at $x'_0$; this shows that
$\phi_0$ is the ``best'' approximation of $\phi$ at $x'_0$ by a projective
map, and therefore that
the third derivative of $\phi_0^{-1}\circ \phi$ at $x'_0$ determines the
Schwarzian derivative $\cS(\phi)$ at $x'_0$ (see e.g. \cite{osgood-stowe},
(1.3)). By construction,
$$ (\phi_0^{-1}\circ \phi)_*I'^*=(\phi_0^{-1})_*(\phi_*I'^*) =
(\phi_0^{-1})_*I^* = \phi_0^*I^*~, $$ 
so 
$\phi_0^{-1}\circ \phi$ is an isometry between $I'^*$ and $\phi_0^*I^*$.
Those two metrics on $\dr_iM'$ are conformal, we will show that the Hessian at 
$x'_0$ of the conformal factor is $\phi^*\II^*_0$, and deduce from this the
proof of Lemma \ref{lm:schwarzian}.

\begin{prop} \label{pr:B}
The shape operator of $\Sigma_\rho$ is given by $B_\rho =
(E+e^{-2\rho}B^*)^{-1} (E-e^{-2\rho}B^*)$, so that, as $\rho\rightarrow
\infty$, $B_\rho - E \simeq -2e^{-2\rho}B^*$. 
\end{prop}

Note that $B_\rho$ and $B^*$ are defined at different points ($B_\rho$ on
$\Sigma_\rho$ and $B^*$ on $\dr_iM$) which are implicitly identified here
through the Gauss flow. 

\begin{proof}
We have already computed above that 
$$ I_\rho=(1/2)(e^{2\rho}I^* + 2\II^* + e^{-2\rho}\III^*) = (1/2)I^*((e^\rho
E + e^{-\rho} B^*)\cdot, (e^\rho E + e^{-\rho} B^*)\cdot)~. $$ 
It follows that
$$ \II_\rho = \frac{1}{2}\frac{dI_\rho}{d\rho} = \frac{1}{2}(e^{2\rho}I^* -
e^{-2\rho}\III^*) = \frac{1}{2} I^*((e^\rho
E + e^{-\rho} B^*)\cdot, (e^\rho E - e^{-\rho} B^*)\cdot)~, $$
and the result is a direct consequence.
\end{proof}

\begin{df}
For each $x$ in some neighborhood of $x_0$ and all $\rho$ large enough, we
call $d_\rho(x)$ the oriented distance, along $\gamma_x$, between its
intersection with $\Phi_0(\Sigma'_\rho)$ and with $\Sigma_\rho$. Then we define
$d_\infty(x):=\lim_{\rho\rightarrow \infty} d_\rho(x)$. 
\end{df}

\begin{prop} \label{pr:distance}
The limit distance $d_\infty$ vanishes at $x_0$ along with its derivative,
and its Hessian at $x_0$ is equal to $\II^*_0$.
\end{prop}

\begin{proof}
By construction, $\gamma_x$ is orthogonal to $\Sigma_\rho$. But
$\Phi_0(\Sigma'_\rho)$ is tangent to $\Sigma_\rho$ at their common intersection
with $\gamma_{x_0}$, so that $d_\rho$ vanishes along with its derivative at
$x_0$. 

Moreover, the second-order variation of $d_\rho$ at $x_0$ is determined by the
difference between the shape operators at their intersection point with
$\gamma_x$ (for $x$ close to $x_0$) of $\Phi_0(\Sigma'_\rho)$
and of $\Sigma_\rho$. But, by Proposition \ref{pr:B} applied to $\Sigma_\rho$
and to $\Sigma'_\rho$, the dominant term in the difference between those
shape operators is
$$ B'_\rho - B_\rho \simeq e^{-2\rho}(2B^*-E)~. $$
Integrating this second derivative of $d_\rho$ shows that the Hessian of
$d_\rho$ is equal to $I_\rho(e^{-2\rho}(2B^*-E)\cdot, \cdot)$, which is
equivalent to $(1/2)(2\II^*-I^*)$, i.e., to $\II^*_0$. 
The result then follows by taking the limit
as $\rho\rightarrow \infty$.
\end{proof}

Since $\phi_0$ is a projective map, it is conformal, and therefore the metrics
$I'^*$ and $\phi_0^*I^*$ on $\dr_iM'$ are conformal, which means that there
exists a function $u:\dr_iM'\rightarrow \R$ such that
$I'^*=e^{2u}\phi_0^*I^*$. 

\begin{prop} \label{pr:factor}
$u(x'_0)=0$ and $du=0$ at $x'_0$, while the Hessian of $u$ at $x'_0$ is equal
to $\II^*_0$.
\end{prop}

\begin{proof}
First note that by definition of $I^*$ it can be described in terms of the
Gauss map, which we call $G_\rho$, from $\Sigma_\rho$ to $\dr_iM$, sending a
point $x\in \Sigma_\rho$ to the endpoint at infinity of the geodesic ray
starting from $x$ orthogonal to $\Sigma_\rho$:
$$ I^* = (1/2)e^{-2\rho} G_{\rho *}(I_\rho +2\II_\rho+\III_\rho) \simeq
2e^{-2\rho} G_*I_\rho~. $$
The same can be said of $\Phi_0(\Sigma_\rho)$, using the Gauss map $G'_\rho$
of $\Sigma'_\rho$:
$$ I'^* \simeq (1/2) e^{-2\rho} G'_{\rho *}(I'_\rho
+2\II'_\rho+\III'_\rho)\simeq 2e^{-2\rho} G'_{\rho *}I'_\rho~. $$ 
Considering the image of $\Sigma'_\rho$ by $\Phi_0$ and the behavior of the
geodesic rays orthogonal to $\Phi_0(\Sigma'_\rho)$ then shows that, in the
neighborhood of $x_0$,
$$ \phi_{0*}(G'_{\rho *}I'_\rho) \simeq e^{2d_\rho}G_{\rho *}I_\rho~, $$
and the result follows by taking the limit as $\rho\rightarrow \infty$.
\end{proof}

\begin{proof}[Proof of Lemma \ref{lm:schwarzian}]
Let $f:=\phi_0^{-1}\circ \phi$. 
By construction the map $f:\dr_iM'\rightarrow \dr_iM'$
sends $x'_0$ to $x'_0$, and it is tangent to the identity at $x'_0$. 
We have
seen above that $f$ is an isometry between $I'^*$ and $\phi_0^*I^*$, so it
follows from Proposition \ref{pr:factor} that $|f'(z)|^2 = e^{-2u}$, where 
$u(x'_0)=0$, $du=0$ at $x'_0$, and the Hessian of $u$ at $x'_0$ is
$\II^*_0$. 

Since $(\dr_iM', I'^*)$ is hyperbolic, we can identify a neighborhood of
$x'_0$ with a neighborhood of $0$ in the unit disk. Since $f$ is conformal it
is then holomorphic, and can be written as: 
$$ f(z) = z + \frac{b}{2}z^2 + \frac{c}{6}z^3 + o(z^3)~, $$
with $b=f''(0)$ and $c=f'''(0)$, and clearly $|f'(z)|=e^{-u}$. But the fact that
$du=0$ at $x'_0$ clearly implies that $b=0$, and it follows that:
$$ |f'(z)|^2 = \left|1+\frac{c}{2}z^2 + o(z^2)\right|^2 
= 1 + Re(cz^2) + o(z^2)~. $$
This means that:
$$ Re(cz^2) = -\II^*_0(z,z)~, $$
where $\II^*_0$ is considered as a real-valued symmetric bilinear form on
$\C$, identified with $T_{x'_0}\dr_iM'$. So
$$ Re((\phi_0^{-1}\circ\phi)'''(x'_0))=-\II^*_0~, $$
which means precisely that $2\II^*_0=-Re(\cS(\phi))$, as claimed.  
\end{proof}


\def\cprime{$'$}

\end{document}